\documentclass[12pt]{article}
\usepackage{MnSymbol}
\usepackage{amsmath,amsfonts,
verbatim}

\newcommand{\Q}{{\mathbb Q}}
\newcommand{\C}{{\mathbb C}}

\newcommand{\PP}{{\mathbf P}}

\newcommand{\LL}{\mathcal{L}}

\newcommand{\F}{\mathbb{F}}
\newcommand{\U}{{\mathbb U}}

\newcommand{\X}{{\mathbb X}}
\newcommand{\XX}{\tilde{\mathbf{X}}}
\newcommand{\Y}{\mathbb{Y}}

\newcommand{\Spec}{\mathrm{Spec}\,}
\newcommand{\M}{\mathbf{M}}   
\newcommand{\MM}{\mathcal{M}}
\newcommand{\NN}{\mathbf{N}} 
\newcommand{\la}{\langle}
\newcommand{\ra}{\rangle}
\newtheorem{pkt}{}[section]  
\newcommand{\bpk}{\begin{pkt}\rm }  
\newcommand{\epk}{\end{pkt}} 
\newcommand{\inv}{^{-1}}   

\newcommand{\be}{\begin{equation}}  
\newcommand{\ee}{\end{equation}}

\newcommand{\dcl}{\mathrm{dcl}}
\newcommand{\acl}{{\rm acl}}

\newcommand{\K}{\mathrm{K}}
\newcommand{\kk}{\mathrm{k}}

\newcommand{\pp}{\mathbf{p}}
\newcommand{\Aut}{\mathrm{Aut}}

\usepackage{graphicx}

\usepackage{tikz}

\title{Definability, interpretations  and
\'etale fundamental groups}
\author{R.Abdolahzadi and B.Zilber\footnote{Supported by the EPSRC program grant ``Symmetries and Correspondences''}}

\begin{document}
\maketitle
\abstract{The aim of the paper and of a wider project is to translate main notions of anabelian geometry into the language of model theory. Here we finish with giving the definition of the \'etale fundamental group $\pi^{et}_1(\X,x)$ of a non-singular quasiprojective scheme  over a field of characteristic 0. }  

\section{Introduction}
Grothendieck's anabelian geometry has been introduced in the language of schemes.
The task of translating its notions into the setting of model theory is not fully trivial for the reason that schemes are essentially objects of a syntax and no semantics is being provided by definitions except some appeal to the intuition of  Zariski-style algebraic geometry. In furnishing a semantics in Tarski's sense one has to take into account a new phenomenon, the syntax here has its own inherent algebraic structure, may have its own automorphisms and homomorphisms, and semantics needs to reflect all those.

 A similar, parallel phenomenon takes place in non-commutative geometry, where syntax is given by ``co-ordinate'' algebras which not only have a structure in algebraic sense but also may be  topological algebras, $C^*$-algebras and so on.   An approach to semantics of non-commutative geometry was suggested by the second author in \cite{Zbook} Ch.5, and \cite{ZQM}. Our approach to schemes is partially inspired by this.  Basic notions are clarified in section \ref{s2}.

In the current paper we introduce a language $\mathcal{L}_\X$ for universal and \'etale covers of a scheme $\X$  and describe the first order theory $T_\X$  of  ``universal \'etale cover'' $\XX^{et}$
of $\X.$ We prove that universal analytic cover $\XX^{an}$ of the complex variety $\X(\C)$ presented in the same language is another model of the theory.  Note that previously various attempts to find an adequate language $\mathcal{L}_\X$ were made in \cite{Zrav}, \cite{Misha}, \cite{DH} and some other publications, covering some classes of varieties $\X.$ 

A key feature of Grothendieck's anabelian geometry is the functor from the Galois category of \'etale covers into the category of groups. In model theory setting this corresponds to the functor $\M\mapsto \Aut(\M)$ from the category of structures with interpretations in the role of morphisms to the category of  topological groups. This functor is well-known in model theory but we give it a more systematic treatment in  section \ref{s1}, especially for the category of finitary  
structures when it becomes an equivalence of categories. In particular, closed subgroups of $\Aut(\M)$ correspond to sets of imaginary elements of $\M.$  

We are especially .interested in classifying imaginaries which correspond to a section of the restriction homomorphism
$\Aut(\M)\to \Aut(\NN)$ when $\NN$ is a structure interpretable in $\M.$ 
 
In the final sections of the paper we give a 
model-theoretic definition of $\pi^{et}_1(\X,x)$ as an automorphism group of $\XX^{et}$ acting  on the fibre over $x.$ We also remark that this group is isomorphic to the Lascar group of the theory. This allows us to reformulate the Grothendieck section conjecture in terms of the properties of structure   $\XX^{et}.$ 

\section{Interpretations and automorphism groups} \label{s1}

\bpk Most of the material in this section is known in some forms. See \cite{CDM} and \cite{MTB} for a model-theoretic approach which we further pursue here. The community of anabelian geometers prefers to speak in terms of Galois categories, see e.g. \cite{C}. One of the aims of the current project is to demonstrate advantages of the model-theoretic point of view. 

Unlike the above publications we do not apriori  restrict the power of the language to first order. The default assumptions is that

\be \label{def=inv} \mbox{\em A relation is definable iff it is invariant under  automorphisms}\ee

For finite structures this property holds for first order languages. For countable structures the language $L_{\omega_1,\omega_1}$ serves the purpose. The main interest to us are {\em finitary} structures defined below. For this class of structures first-order languages are essentially sufficient. 

{\bf Definable} means definable without parameters (the same as 0-definable). 

In general, we consider multi-sorted $\LL$-structures $\M$. A definable set in 
$\M$ is a definable subset $D$ of $S_1\times \ldots \times S_n,$ a cartesian product  of finitely many sorts.

A definable {\bf sort} in an $\LL$-structure $\M$ is a set of the form $D/E$ where $D=D(\M)$ is a definable set in $\M$ and $E$ a definable equivalence relation on $D.$ An $n$-ary relation  on  $D/E$ is definable if its pull-back under the canonical map
$D\to D/E$ is  definable. 

{\bf An interpretation} of an  $\LL_N$-structure $\NN$ in an $\LL_M$-structure $\M$ is a bijection
$g: N\to D/E,$  a sort in $\M$ such that 
for each basic relation $R$ the image $g(R)$ is a definable relation on the sort  $D/E.$

Given a structure $\M$ we also consider the structure $\M^{Eq}$ interpretable in $\M$ and which has every sort of $\M$ as a definable substructure.

Note that a finite union of sorts and  a direct product of finitely many sorts is a sort in $\M^{Eq}.$
 
\epk
\bpk {\bf Category $\mathfrak{M} .$} Its objects are (multisorted) $\LL$-structures $\M$ (all $\LL$).
 
 The {\bf pre-morphisms} $g:\mathbf{N}\to \M$
 are interpretations (without parameters).  More precisely,
 $$g: \mathbf{N}\to \M^{Eq}$$ is an injective map such that $g(N)$  is a universe of a sort in $\M^{Eq},$  and for  any basic relation or operation $R$ on  $\mathbf{N}$ the image $g(R)$ is definable  in the sort. 
 
 We denote $g(\mathbf{N})$ the  $g(N)$ together with all the relations and operations $g(R)$ for $R$ on $\NN.$

\medskip
 
Two pre-morphisms $g_1: \NN\to \M$ and $g_2\to \NN\to \M$ are equivalent if  there is a bijection $h:g_1(\NN)\to g_2(\NN)$   which is definable in $\M^{Eq}.$

The equivalence class of a pre-morphism $g:\NN\to \M$ is a {\bf morphism} $\mathbf{g}:\NN\to \M.$  

\medskip

The following definitions will be used  for pre-morphisms $g$  as well as for morphisms
$\mathbf{g}.$
 
 We say $g$ is an {\bf embedding},  $g:\mathbf{N}\hookrightarrow \M$ if $g(\mathbf{N})$ has no proper expansion definable in $\M^{Eq}.$

We say $g$ is a {\bf surjection}, $g:\NN\twoheadrightarrow \M$ if $M\subseteq \dcl(g(N))$
where $\dcl$ is in the sense of $\M^{Eq}.$

We say that $g: \NN\to \M$ is an {\bf isomorphism},  $g:\mathbf{N}\cong \M,$ if
$g$ is an embedding and a surjection.
 
\medskip

In what follows we sometimes write $\mathbf{N}\cong_{\mathfrak{M}} \M$ to emphasise that the isomorphism (or morphism) is in the sense of the category $\mathfrak{M}$ to distinguish from ones in the usual algebraic sense.

\medskip

{\bf Lemma.} {\em Let $g:\mathbf{N}\to \M$ be an $\mathfrak{M}$-isomorphism and let $\M'=g(\NN).$ Then the inverse map $g^{-1}:\M'\to \mathbf{N}$ induces 
 a $\mathfrak{M}$-isomorphism 
$h: \M\to \NN.$}
\medskip

{\bf Proof.} 
By assumptions
we have  $M\subseteq \dcl (M')$ in $\M^{Eq}.$ This implies that there are in $\M:$
a family $\{ S_i: i\in I\}$   of definable
subsets $S_i\subset  {M'}^{n_i}$ and 
 a family of definable functions $h_i: S_i\to M$ such that $$\bigcup_{i\in I} h_i(S_i)=M\mbox{ and }h_i(S_i)\cap h_j(S_j)=\emptyset \mbox{ if } i\neq j .$$

{\bf Claim.} We may assume that the family $\{ S_i: i\in I\}$ of domains of $f_i$ is disjoint, that is $S_i\cap S_j=\emptyset$ if $f_i\neq f_j.$

Proof. Note that by definition $\dcl(M')=\dcl(M'\cup \dcl(\emptyset)),$ where $\dcl$ is understood in the sense of $\M^{Eq}.$ The latter has, for each $i\in I,$ the sort '$f_i$' which is defined as the $\mathrm{graph} f_i/E_i$ where $E_i$ is the trivial equivalence relation with one equivalence class. Clearly, '$f_i$'$\in \dcl(\emptyset).$ 
 Now replace 
$S_i$ by $S_i\times$'$f_i$' and we have the required.
   
   \medskip

Set $D(\M'):= \bigcup_{i\in I} S_i$ and 
$h: D(\M')\to \M$ to be  $\bigcup_{i\in I}h_i.$
 $h$ is a map definable in $\M$ and is an interpretation, a pre-morphism $\M\to \M'.$ On the other hand, any relation  on $\M'$ is a relation on a sort in $\M$ since $\M'$ is a sort in $\M^{Eq},$   hence there are no new relations on $h(M)$ induced from $\M',$ that is the interpretation $h$ is an embedding.  Recalling that $\M'= g(\NN)$ completes the proof. 
 $\Box$ 

\medskip

We identify morphism $h$ as in the Lemma with $g\inv.$
 \epk
\bpk For a subset $A\subseteq M,$ denote $\M/A$ (or sometimes $\M(A)$) the expansion of $\M$ by  names of elements of $A.$

Clearly, the identity map defines a (canonical) morphism $\M\to \M(A).$ This morphism is an embedding (and so isomorphism) if and only if $A\subseteq \dcl(\emptyset).$
\epk
\bpk \label{dcl} Given $A\subseteq \dcl(\emptyset)$ we may treat $A$ as a structure in which any element is named (e.g. by a formula defining the element in $\M$) and so any relation is definable. Clearly then
$$\mathrm{Aut}(A)=1 \mbox{ and } A\hookrightarrow \M.$$
 
\epk

 \bpk The category  $\mathfrak{M}_\mathrm{fin}$ is a subcategory of  $\mathfrak{M} $ whose objects are {\bf finitary} structures $\M,$ that is structures which can be represented in the form
 $$\M=\bigcup_{\alpha<\kappa } \M_\alpha$$
 where the $\M_\alpha$ are finite first-order $0$-definable substructures of $\M.$
 \medskip

 Note that an equivalent definition would be 
 $$\M=\acl(\emptyset)$$
 where $\acl$ is in the sense of first-order logic.
\medskip

 {\bf Example.} Let $\kk$ be a field and $\F=\tilde{\kk},$ its algebraic closure. We consider $\F=\F/\kk$ as a structure
 in the language of rings with names for elements of $\kk.$ Then each $a\in \F$ is contained in a $0$-definable set $M_a$ equal to its Galois orbit $M_a:=G_k\cdot a,$ $G_k=\mathrm{Gal}(\F:\kk).$  So  $\F/\kk\in \mathfrak{M}_\mathrm{fin}.$
 \epk
 
 \bpk \label{EqCat}{\bf Theorem.} {\em The map $\M\to \mathrm{Aut}(\M)$ induces a contravariant  functor from
 $\mathfrak{M}$ into the category  $\mathfrak{G}_{\mathrm{top}}$ of topological groups. This functor sends $\mathfrak{M}_\mathrm{fin}$ into the category of profinite groups $\mathfrak{G}_{\mathrm{pro}},$.  
 
(i) To every $g:\mathbf{N}\to \M$ corresponds the restriction homomorphism 
 
\noindent $\hat{g}:\mathrm{Aut}(\M^{Eq})\to \mathrm{Aut}(\mathbf{N}).$
 
 (ii) An embedding $g:\mathbf{N}\hookrightarrow \M$ to the surjection $\hat{g}:\mathrm{Aut}(\M)\twoheadrightarrow \mathrm{Aut}(\mathbf{N}),$ provided $\M\in \mathfrak{M}_\mathrm{fin}.$
 
 (iii) The expansion by naming all points in $A\subseteq \M^{Eq},$
 $g:\M\to \M/A$ corresponds to an embedding $\hat{g}:\mathrm{Aut}(\M/A)\hookrightarrow \mathrm{Aut}(\M).$

(iv) The profinite restriction of    the functor, $$\mathrm{Aut}:\ \mathfrak{M}_\mathrm{fin}\to \mathfrak{G}_{\mathrm{pro}},$$ 
 is an equivalence of categories. 
 

   
}
\medskip

{\bf Proof.} (i)
is immediate by definition.

(ii) Since $g$ is an embedding, the relations definable on $g(\mathrm{N})$ are the same in $\M$ and $\mathbf{N}.$ Hence a $g(\mathbf{N})$- automorphism $\alpha$ is an elementary
bijection $g(\mathrm{N})\to g(\mathrm{N})$ in $\M.$ Now consider an elementary saturated enough extension $\M\prec {^*\M}.$ It is homogeneous and so $\alpha$ extends to an automorphism $\alpha^*$ of  ${^*\M}.$ But any automorphism preserves $\mathrm{acl}_{\M}(\emptyset)$ which is equal to $\M$ under our assumption. Thus $\alpha^*$ induces an automorphism of $\M$ which extends $\alpha.$

(iii) Immediate.

(iv) 
First we prove the statement for $\Aut:\mathfrak{M}_\mathrm{finite}\to \mathfrak{G}_{\mathrm{finite}},$ the functor between finite structures and finite groups, 
subcategories of  $\mathfrak{M}_\mathrm{fin}$ and  $\mathfrak{G}_{\mathrm{pro}},$ respectively.

Given a finite group $\mathbf{G}$ one constructs a finite $\M$ such that $\mathbf{G}\cong \mathrm{Aut}(\M)$ by setting $M=G$ and introducing all relations $R$ on $M$ which are invariant under the action of $\mathbf{G}$ on $G$ by multiplication. This gives us for $\M=(M; R)$ 

Claim.
$$\mathbf{G}=\mathrm{Aut}(\M)$$
Proof. $\mathbf{G}$ acts on $\M$ by automorphisms by definition. We need to prove the inverse, i.e. that there are no other automorphisms. Consider the tuple $\bar{g}$ of all the elements of $G$ (of length $n=|G|$) and let $S_g$ be the conjunction of all the relation in $R$ that hold on $\bar{g}$ (that is $\mathrm{tp}(\bar{g})$). We can also consider $S_g^0:=\mathbf{G}\cdot \bar{g},$ the orbit of $\bar{g}$ under the action of $\mathbf{G}.$ Clearly,
$S_g^0\subseteq S_g$ and by minimality they are equal. 

Now take an automorphism $\sigma$ and consider $\sigma \bar{g}.$ This is in $S_g$ and thus, for some $h\in G,$    $\sigma \bar{g}=h\bar{g},$ that is $\sigma g_i=hg_i$ for each $g_i\in G.$
Claim proved.

It remains to see that if  $\mathbf{G}\cong\mathrm{Aut}(\mathbf{N}),$ then   $\mathbf{N}$ is definable in $\M$ and vice versa. In order to do this we may assume $\mathbf{G}=\mathrm{Aut}(\mathbf{N}).$  

Consider $N,$ the universe of the structure, and let $\mathbf{n}$ be the $N$ presented as an ordered tuple. 
Let $M':= \mathbf{G}\cdot \mathbf{n},$ the orbit of the tuple under the action of the automorphism group. Clearly, $M'$ consists of $|\mathbf{G}|$ distinct elements, since automorphisms differ if and only if they act differently on the domain $N.$ Also $M'$ is definable in $\mathbf{N}$ since the tuples $\mathbf{n}'$ making up $M'$ are characterised by the condition that $\mathrm{tp}(\mathbf{n}')=\mathrm{tp}(\mathbf{n}).$ The relations $R$ induced on $M'$ from $\mathbf{N}$ are invariant under $\mathrm{Aut}(\mathbf{N}),$ and because a finite structure is homogeneous, the converse holds. In other words an obvious bijection $M\to M'$ is a bi-interpretation, so $\M\cong \M'$ in the sense of $\mathfrak{M}.$  At last notice that we can interpret $\mathbf{N}$ in $\M'$ since the relation  ``$\mathbf{n}'$ and $\mathbf{n}''$ have the same first coordinate is invariant under $\mathbf{G}$'' is definable. This gives us $N$ as a definable sort. It follows that any relation on $N$ definable in $\mathbf{N}$ is definable in $\M'.$ So $\mathbf{N}\cong \M'\cong \M$  in the sense of $\mathfrak{M}.$  Finite case of $\Aut$ proven.

Now we extend $\Aut$ to the category of  finitary structures   $\M\in \mathfrak{M}_{\mathrm{fin}}$ by continuity
$$\M=\lim_\to \M_\alpha \to \mathbf{G}=\lim_{\leftarrow}  \mathbf{G}_\alpha, \mbox{ where } \mathbf{G}_\alpha=\Aut  \M_\alpha$$

Since the functor is invertible and preserves morphisms on finite objects of the categories, it is an equivalence also on the limits.

$\Box$
 \epk
\bpk {\bf Example.} Let $K$ and $L$ be two number fields, $\tilde{\Q}=\F.$ Let $\F_K$ and $\F_L$ be two structures with respective subfields of constants (named points). Clearly 
these belong to $\mathfrak{M}_\mathrm{fin}.$
A celebrated theorem by  Neukirch states that
$$\F_K\cong_\mathfrak{M} \F_L\Leftrightarrow K\cong L.$$

\epk
\bpk \label{[N]} Let  $i: \mathbf{N}\to \M/A$ be an interpretation of $\NN$ in $\M$ over $A.$ Then every relation $R$ which is $0$-definable in  $\mathbf{N}$  is a point in $\M^{Eq}/A.$  The definable closure of all these points  denote $[i\mathbf{N}]$ or often just   $[\mathbf{N}]$ where an $i$ is assumed.

More precisely, let

\be\label{N2} \mathrm{Aut}(\M/ [i\mathbf{N}]):=\{ \sigma\in \mathrm{Aut}(\M^{Eq}): \sigma_{|i\mathbf{N}}\in \mathrm{Aut}(i\mathbf{N})\},\ee

Set
$$[i\NN]:= \mathrm{Fix}_{\M^{Eq}}(\mathrm{Aut}(\M/ [i\mathbf{N}])=\{ a\in \M^{Eq}: \forall \sigma\in \mathrm{Aut}(\M/ [i\mathbf{N}])\ \sigma(a)=a \} .$$

Clearly, \be\label{N1} [i\mathbf{N}]\subseteq \dcl(A).\ee

\medskip

Also,

and if    $i\mathbf{N}$ is $0$-definable then 
$$[i\mathbf{N}]=\dcl(\emptyset)\mbox{ and }\mathrm{Aut}(\M/[i\mathbf{N}])=\mathrm{Aut}( \M).$$
\epk

 
\bpk\label{sections}
 {\bf Proposition.} {\em Let $\mathbf{N}, \M\in \mathfrak{M}_{\mathrm{fin}}.$
 Let $\hat{h}:\mathrm{Aut}(\M)\twoheadrightarrow\mathrm{Aut}(\mathbf{N})$ and let 
$\hat{g}:\mathrm{Aut}(\mathbf{N})\hookrightarrow \mathrm{Aut}(\M)$ be a section of 
$\hat{h},$ that is 
$\hat{h}\circ \hat{g}=\mathrm{id}_{\mathrm{Aut}(\mathbf{N})}.$

Then there is an interpretation-isomorphism $$i: \mathbf{N}\cong \M/A,$$
for some $A\subset M$ satisfying
$\dcl_{\M^{Eq}}(A)=\dcl_{\M^{Eq}}([i\NN]),$ so we may set $A=[i\NN],$ such that
 \be\label{sections1} M\subseteq\dcl_{\M^{Eq}}(iN\cup [i\NN])  \mbox{ and } [i\NN]\cap iN =\dcl_{\mathbf{N}}(\emptyset).\ee
   }
   
   {\bf Proof.} By \ref{EqCat} and the assumptions, there is an interpretation 
  $$h: \NN\hookrightarrow \M$$ 
  correponding to $\hat{h},$ such that any $\sigma\in \Aut(\M)$ induces  $\hat{h}(\sigma)\in \Aut(h\NN)$ and in this way we get all automorphisms of $h\NN,$ that is $\hat{h}(\Aut(\M))=\Aut(h\NN).$ Without loss of generality we may assume that $\NN$ is a substructure of $\M^{Eq},$ that is
  $h$ is a pointwise identity embedding and $\hat{h}(\sigma)$ is the restriction of $\sigma$ to $\NN.$ Thus 
  
  \be\label{hN} \hat{h}(\Aut(\M))=\Aut(\NN).\ee 
    
   Consider the subgroup $\hat{g}(\Aut(\NN))\subseteq \Aut(\M),$ an isomorphic copy of $\Aut(\NN).$ Since $\hat{g}$ is a section of $\hat{h}$ we get $$\hat{h} (\hat{g}(\Aut(\NN)))=\Aut(\NN).$$
   
By assumptions $\hat{g}$ lifts any automorphism
$\rho\in \Aut(\NN)$ to a unique automorphism
$\hat{g}(\rho)\in \Aut(\M)$ giving the embedding $\hat{g}:\Aut(\NN)\hookrightarrow \Aut(\M).$ 
   
   
Set $$A:= \mathrm{Fix}_{\M^{Eq}}(\hat{g}(\Aut(\NN))).$$

Note that according to  (\ref{N2}) 
$$\hat{g}(\Aut(\NN))=\Aut(\M/[\NN])$$  
   and this by definition in \ref{[N]} 
   $A=[\NN].$  Moreover, $\NN$ is definable in $\M^{Eq}$ over $A$ since $\NN$ as a structure is  $\Aut(\M^{Eq}/A)$-invariant. In other words there is a pre-morphism $i:\NN\to \M.$ 
   
Now note that $i$ is an $\mathfrak{M}$-embedding since every  $\rho\in\Aut(\NN)$ lifts to an automorphism $\hat{g}(\rho)$ of $\M^{Eq}/A.$

Next we note that $i$ is an $\mathfrak{M}$-surjection, that is $\dcl_{\M^{Eq}/A}(N)\supseteq M,$ or equivalently  $$\dcl_{\M^{Eq}}(N\cup A)\supseteq M.$$ To see the latter we remark that if $\sigma\in \Aut(\M)$ fixes $A\cup N$ point-wise then $\sigma\in \hat{g}(\Aut(\NN))$ (because $A$ is fixed) and $\sigma$ is identity on $\NN.$ That is $\sigma=1.$  

Finally note that $A\cap N=\dcl_{\NN}(\emptyset)$        
because the intersection consists exactly of $\Aut(\NN)$- fix-points of $\NN.$ $\Box$  
\epk

\bpk \label{M/N} {\bf Lemma.} {\em Suppose $\M$ is algebraic over $\emptyset,$ $\M=\acl(\emptyset).$
Let $G\hookrightarrow \mathrm{Aut}(\M)$ be a closed subgroup (in the profinite topology). Then $G$ is a pointwise stabiliser of a subset $A\subset M.$  That is $$G=\mathrm{Aut}(\M/A)$$

$G$ is normal iff $A$ can be chosen $0$-definable.}

\medskip

{\bf Proof.} Let $A=\{ a\in M: \ \forall g\in G\, g.a=a\}.$ Then $G\subseteq \mathrm{Aut}(\M/A).$ The inverse follows from $G$ being closed. This is easy to see for $\M$ finite, and this is enough
 since $\M$ is algebraic.

With the above choice of $A,$  $G$ is normal iff $N$ is invariant under  $\mathrm{Aut}(\M),$ that is
$L_{\omega_1,\omega}$-definable. $\Box$
\epk

\bpk\label{Pr1} {\bf Proposition.} {\em To every $0$-definable  $\NN$ in $\M$ (write $\NN\hookrightarrow \M$)
one associates the exact sequence
$$1\to \mathrm{Aut}(\M/\NN)\to \mathrm{Aut}(\M)\to \mathrm{Aut}(\NN)\to 1,$$
and every exact sequence of closed subgroups has this form for some  $\NN\hookrightarrow \M.$}
\epk

 \section{ Reduced schemes and varieties over $\kk$}\label{s2}

Our aim in this section is to clarify the relations between scheme-theoretic language and the language of varieties which is more readily adaptable to model theory treatment. 

By  scheme we always mean a reduced scheme.

Recall that an affine $\kk$-scheme of finite type is given by a  commutative unitary ring $A$  with $\kk\hookrightarrow A,$ finitely generated over $\kk$ and without nilpotent elements. From logical perspective this should be treated as a syntax to which we still have to provide semantics. One can do it by first, associating with every $A$ a certain language $L^A$, and then providing a first-order $L^A$-theory $T^A,$  models of which will be seen as semantic realisations of $A,$ not necessarily the scheme as $\Spec A$ in scheme-theoretic sense.    The morphisms $A\to B$ (in particular, automorphisms) given by homomorphisms of rings in the category of  schemes must be reflected by certain ``morphisms'' between models of $T^A$ and $T^B.$

\bpk  \label{defV} {\bf Varieties over $\kk.$}
Let $\F$ be an algebraically closed field containing $\kk$ and $\F_{/\kk}$ be the field $\F$ with names for all elements of $\kk.$

A {\bf (formal) affine variety $\X$ over $\kk$ realised in $\F$} (sometimes written as $\X_{/\kk}(\F)$ or $\X_{\kk}$ )
is the two sorted structure
$\left(\X(\F), \F_{/\kk}\right),$ in a language $L^A,$ where $\X(\F)\subseteq \F^n$ is the zero set of a system of polynomials over $\kk$ (the variety sort), $\F$ is the sort for the field, and $L^A$ is the language with
\begin{itemize}
\item  unary predicates for sorts $\X$ and $\F$,
\item  the addition and multiplication operation on sort $\F,$
\item constant symbols for each element of $\kk$ in $\F,$
\item names of all the Zariski regular unary maps $a: \X(\F)\to \F$ defined over $\kk.$\footnote{If $\X=\X_1{\cup}\ldots {\cup} \X_l$ decomposition into irreducible components, $a: \X(\F)\to \F$ is regular if and only if its restrictions on each component is regular.} 
\end{itemize}

 Note that since $\X(\F)\subseteq \F^n$ the genuine $n$ coordinate functions, projections into $\F,$ generate the $\kk$ algebra of regular functions $\X(\F)\to \F$ conventionally written $\kk[\X].$

Fixing $x\in \X(\F)$ one gets a homomorphism (representation) of the $\kk$-algebra $A$ into the $\kk$-algebra $\F_{/\kk}$
\be \label{rho0}\rho_x: a\mapsto a(x),\ \ A\to \F,\ee
so the points of $\X(\F)$ are in a natural bijective correspondence with irreducible representations of $A$ and  in a natural bijective correspondence with maximal ideals of $A.$

Choosing a finite set $a_1,\ldots,a_m\in A$ generating $A$ as an algebra over $\kk$ one gets an emdedding
$$x\mapsto \la a_1(x),\ldots,a_m(x)\ra;\ \ \X(\F)\hookrightarrow \F^m$$
\epk
\bpk To every connected (not necessarily geometrically connected) affine $\kk$-scheme $A$ we now associate a first-order $L^A$-theory $T^A.$ 

Let $\kk^A\subset A$ be a maximal subfield of $A$ algebraic over $\kk.$ Since $A$ is finitely generated, $\kk^A=\kk[a]$ for some $a\in A.$

Note that there is just one maximal subfield of $A$ algebraic over $\kk.$ 
This follows from the fact that $A$ has no zero-divisors.

\medskip

The theory $T^A$ of two sorted structures $(\X,\F)$
will say:
\begin{itemize}
\item[1.] $(\F,+,\cdot)$ is an algebraically closed field;
\item[2.] For each $\xi\in \kk\subseteq \F,$ there is $a\in A$ such that 
$ \forall x\in \X \  a(x)=\xi$;

For each $a\in \kk^A$ there is $\alpha\in \F$ such that $ \forall x\in \X \  a(x)=\alpha;$

\item[3.] For each non-constant polynomial $p(T_1,\ldots,T_m)$ over $\kk,$ for any $a_1,\ldots,a_m\in A,$   $$p(a_1,\ldots,a_m)=0\mbox{ iff }\forall x\in \X\ p(a_1(x),\ldots,a_m(x))=0;$$
\item[4.]  
For each  $a_1,\ldots,a_n\in A\otimes_{\kk^A=\kk[\alpha]} \F,$ such that the ideal $\la a_1,\ldots,a_n\ra$ 
does not contain 1, it holds $$\exists x\in \X\   a_1(x)=\ldots =a_n(x)=0.$$

\item[5.] Given a  generating set $a_1,\ldots,a_n$ of $A,$
$$\forall x_1,x_2\in \X\ (\bigwedge_{i=1}^n a_i(x_1)=a_i(x_2) \ \to \ x_1=x_2).$$
\end{itemize}

Axioms 2. - 4. are infinite sets of first-order axioms each. 
\epk
\bpk \label{schemes1} {\bf Proposition.} {\em $T^A$ is a complete theory categorical in uncountable cardinals. Any model of $T^A$  is isomorphic to one of the form $(\X(\F),\F_{/\kk})$ for some affine variety $\X(\F)\subseteq \F^n$ such that the co-ordinate algebra $\kk[\X]$ is isomorphic to $A.$

}

{\bf Proof.} Given $T^A,$ axioms 2 and 3 determine $A$ as an algebra over $\kk.$

Choose an algebraically closed field $\F$ and choose an embedding $\imath:\kk\hookrightarrow \F.$ 
Interpret of $\imath(\xi)$ as the constant map with value $\xi$  according to the first part of axiom 2. 

Let $\kk^A:=\kk[a].$ Choose an embedding $\imath_a:\kk[a]\hookrightarrow \F$ extending $\imath.$ Let $\alpha=\imath_a(a).$  Interpret $a$ as the constant map with value $\alpha.$

Set \be\label{X(F)} \X(\F):= \mathrm{Hom}_\F(A\otimes_{\kk^A}\F ,\F),\ee where $\kk^A$ is identified with 
$\imath_a(\kk^A)\subset \F.$
 
Elements  $a\in A$ interpret a functions $a: \X(\F)\to \F$ by defining, for $x\in \X(\F),$
$$a(x):=x(a).$$ 
This gives us a model of $T^A.$

Now we note that every model of $T^A$ is of this form. Indeed, by axiom 1 we start with an algebraically closed $\F$ for the respective sort.
Axioms 2 and 3 tell us that algebraic relations in the $\kk$-algebra $A$ agree with
algebraic relations between the  functions $X\to \F$ interpreting elements of $A.$ 
 Axiom 4 guarantees that every 
homomorphism $A\to \F$ gives rise to a point in $\X.$ And axiom 5 implies that the correspondence between points of $\X$ and homomorphisms is bijective. 

It remains to prove categoricity of $T^A$ in uncountable cardinals. We prove a stronger statement, that  any two models with the same $\F$ are isomorphic. Indeed, once $\kk^A=\kk[a]$ is identified with $\kk[\alpha]\subset \F,$ and $a$ with $\alpha,$ the $\F$-algebra $A\otimes_{\kk^A}\F$ is a  coordinate algebra of a unique connected variety $\X(\F)$ over $\F.$ We write it as $\X_\alpha(\F).$  Any other model with the same $\F$ is obtained by choosing an embedding  $\kk[a]\to \F,$   that is by identifying $a$ with an $\alpha'\in \F,$ Galois conjugated to $\alpha.$ By construction (\ref{X(F)}), a field isomorphism $\F_{/\kk[\alpha]}\to  \F_{/\kk[\alpha']},$ $\alpha\mapsto \alpha',$
induces an isomorphism 
$\X_\alpha(\F)\to \X_{\alpha'}(\F)$ which gives us an isomorphism between the two models.
$\Box$

\medskip

Note that given $T^A$ we are given the language $L^A$ and the axioms 2. and 3. provide the description of  the ring structure on $A$ as well as embedding $\kk\to A.$ Hence $T^A$ uniquely recovers $A.$

This describes a bijective correspondence between the class of affine geometrically connected reduced $\kk$-schemes and theories of the form $T^A,$ equivalently, between schemes and elementary classes of  $L^A$-structures of the form $(\X(\F), \F_{/\kk}),$ or more precisely 
 $(\X_\alpha(\F), \F_{/\kk}).$

\epk

\bpk \label{k-var} {\bf Remark.} By associating $\alpha\in \F$ with $a\in A$ we  identify 
$\X_\alpha(\F)$ 
as a connected affine variety defined over $\kk_\alpha:=\kk[\alpha],$ a $\kk_\alpha$-variety. 

The respective embedding $\kk[\alpha]\hookrightarrow A$ allows us to see $\X_\alpha(\F)$ as a $\kk[\alpha]$-scheme, a connected component of scheme $A.$

\epk
\bpk \label{morphisms} Morphisms $f: A\to B$ between affine $\kk$-schemes are homomorphisms of $\kk$-algebras.

 $f$ canonically  induces an interpretation  of models $(\X_\alpha(\F), \F_{/\kk})$ of $T^A$ in appropriate models $(\Y_\beta(\F), \F_{/\kk})$ of $T^B.$ 
 
 Let $\kk^A=\kk[a]$ be as above and
$\kk^B=\kk[b]$ be the algebraic closure of $\kk$ in $B.$

We will have $f_{|\kk[a]}$ an embedding into $\kk[b]$ and hence $f(a)$  is a constant function on  $\Y_\beta(\F),$ $f(a)(y)= \alpha',$ where $\alpha'\in \kk[\beta],$  $\kk[\alpha']\cong \kk[\alpha].$ 

Note that the choice of the constant functions identifies $\kk[\alpha]$
with $\kk^A$ and $\kk[\beta]$
with $\kk^B.$

We will say that the choice of model $(\Y_\beta(\F), \F_{/\kk})$ of $T^B$ is appropriate to the choice of model $(\X_\alpha(\F), \F_{/\kk})$ of $T^A$ with respect to $f$ if $\alpha'=\alpha.$

When $\alpha'=\alpha$ we can define the unique homomorphism of $\F$ algebras
$$f_{\alpha,\beta}: A\otimes_{\kk[\alpha]} \F\to B\otimes_{\kk[\beta]}\F$$
which is given as $f$ on $A$ and as identity on $\F.$

By commutative algebra, there is a canonical regular map between varieties
\be\label{decompf} f^*_{\beta,\alpha}: \Y_\beta(\F)\to \X_\alpha(\F)\ee functorially correponding to $f_{\alpha,\beta}.$

The family of all the maps $f^*_{\beta,\alpha}$ between variety sorts of appropriate pairs of models $T^B$ and $T^A$ we define as the morphism \be\label{f*} f^*:T^A\to T^B\ee
which corresponds to $f:A\to B$ (note that $\to$ in (\ref{decompf}) is inverse to that of (\ref{f*})).

Finally we note that $f^*_{\beta,\alpha}:\Y_\beta(\F)\to \X_\alpha(\F),$ for any choice of appropriate $\alpha,\beta,$ determines $f:A\to B.$ Indeed, given $a\in A$ together with its interpretation $a: \X_\alpha(\F)\to \F$ in the models
we can define $b=f(a)$ uniquely by setting,  
$$b(y)=f(a)(y)= a(f^*_{\beta,\alpha}(y)).$$

\epk

\bpk
In case  $f$ is an embedding, $f^*_{\beta,\alpha}$ is a surjective map and can be  classified as an interpretation  of model $(\X_\alpha(\F), \F)$ of $T^A$ in model $(\Y_\beta(\F), \F)$ of $T^B.$ 
Indeed 
$\X_\alpha(\F)$ can be identified with a quotient of $\Y_\beta(\F)$ by a Zariski closed equivalence relation definable in $L^B.$ Co-ordinate functions $a\in A$ on $\X_\alpha(\F)$ will be interpreted by $f(a)=b\in B$ and hence $L^A$-relations are interpretable in the language $L^B.$

We will say in this case that {\em theory $T^A$ is interpreted in theory $T^B.$}

\medskip

In case $f$ is an embedding and $f(\kk^A)=\kk^B$ for appropriate pairs of models  we have $\alpha=\beta$ and thus both  $\X_\alpha(\F)$ and $\Y_\alpha(\F)$ are $\kk[\alpha]$-varieties. Now the interpretation
of $\X_\alpha(\F)$ in 
$\Y_\alpha(\F)$ does not induce new relations on $\X_\alpha(\F)$ since any automorphism of $\F_{/\kk[\alpha]}$ can be extended uniquely to  (abstract) automorphisms of $\X_\alpha(\F)$ and  $\Y_\alpha(\F)$ respecting $f^*_{\alpha,\alpha}.$

Thus we have proved 

\epk
\bpk\label{schemesP} {\bf Proposition.} {\em The category of affine $\kk$-schemes with morphisms
$f: A\to B$ is equivalent to  the category of theories of varieties with morphisms
$f^*:T^A\to T^B$ described above. The morphisms $f^*$ acts between models as a regular map $f^*_{\beta,\alpha}: \Y_\beta(\F)\to \X_\alpha(\F)$  where $\Y_\beta$ is a $\kk[\beta]$-variety,  $\X_\alpha$ is a $\kk[\alpha]$-variety and there is an embedding $\kk[\alpha]\subseteq \kk[\beta]$ induced by the embedding $f:\kk^A\to \kk^B.$

If $f$ is an embedding, $f^*$ determines an interpretation $f^*:T^A\to T^B$ of appropriate models of the theory $T^A$ in models of the theory $T^B.$

If, in addition, $f(\kk^A)=\kk^B$ the interpretation $f^*$ is an embedding-interpretation.
} 
\epk
The more general notion of a quasi-projective $\kk$-schemes $A$ of finite type is based on the {\bf gluing construction} between several affine schemes $A_1,\ldots,A_n.$ The gluing construction  between schemes can be routinely translated into 
a definable gluing construction between models of $T^{A_1},\ldots,T^{A_n}$ thus furnishing a theory
$T^A$ with two-sorted models of the form $\X(\F), \F_{/\kk},$ where $\X(\F)$ is a quasi-projective variety. Morphisms between models of respective theories are defined accordingly.

We skip the detail and claim the following generalisation of \ref{schemesP} 
\bpk \label{schemesPP} {\bf Theorem.}{\em The category of quasi-projective $\kk$-schemes with morphisms
$f: A\to B$ is equivalent to  the category of theories of varieties with morphisms
$f^*:T^A\to T^B$ described above. The morphisms $f^*$ acts between models as as a regular map $f^*_{\beta,\alpha}: \Y_\beta(\F)\to \X_\alpha(\F)$ where $\Y_\beta$ is a $\kk[\beta]$-variety,  $\X_\alpha$ is a $\kk[\alpha]$-variety and there is an embedding $\kk[\alpha]\subseteq \kk[\beta]$ induced by the embedding $f:\kk^A\to \kk^B.$

If $f$ is an embedding, 
$f^*$ determines an interpretation  $f^*:T^A\to T^B$ of appropriate models of the theory $T^A$ in models of the theory $T^B.$  

If, in addition, $f(\kk^A)=\kk^B$ the interpretation $f^*$ is an embedding-interpretation.} 
\epk

\bpk Let $\X_\alpha(\F)$ be given by its 
embedding $\X_{\mu,\alpha}(\F)\subset \PP^N(\F).$ 
This defines the coordinate algebra over $\F$ of affine charts of the variety and so defines the respective $\F$-scheme, up to isomorphism of schemes.  In particular, 
the embedding defines the respective $\kk(\alpha)$-scheme up to isomorphism of schemes.

Equivalently, defines $T^A$ up to an isomorphism of the language $L^A.$ 
\epk

\bpk \label{rem2.2} {\bf Remark.} (i) Note that for $\X$ over $\kk$
$$\X_{/\kk}(\F)\cong_{\mathfrak{M}} \F_{/\kk}.$$
This is immediate from definitions.


(ii) Besides the obvious interpretation of $\F$ in $\X(\F)$ by the sort $\F$ we have interpretations associated with 
each non-constant coordinate function $a: \X\to \F.$  The field structure $\F$ is interpreted then on the set $\X(\F)/E_a$ (plus-minus a finite subset)
where $E_a$ is an equivalence relation given by $x\, E_a\, x'\Leftrightarrow a(x)=a(x').$  

This interpretation is isomorphic to (bi-interpretable with) the obvious interpretation on the sort $\F$ since the map $a$ induces a definable isomorphism between $\X(\F)/E_a$ and $\F$ and in characteristic 0 the only definable automorphism of $\F$ is the identity.   

\medskip


\epk



\bpk \label{isoSchemes}  {\bf Remark.} Let $\X_{\alpha}(\F)$ and $\X_{\alpha'}(\F)$ be the variety sorts in  two models of $T^A.$ Let  as in \ref{morphisms} $\kk^A=\kk[a]$ where $a$ interpreted as the constant function with value $\alpha$ in $\X_{\alpha}(\F)$ and the constant function with value $\alpha'$ in $\X_{\alpha'}(\F).$ 

Let $f: A\to A$ be an automorphism  such that $f(a)=a'\in \kk^A,$ where $a'$ is the name corresponding to the constant function with value $\alpha'$ in $\X_{\alpha}(\F).$ Then the automorphism $f$ of schemes induces the isomorphism of $f^*_{\alpha,\alpha'}: \X_{\alpha}(\F)\to \X_{\alpha'}(\F)$ of models, and conversely, every isomorphism of models by changing the embedding $\kk^A\to \F$ is induced by an automorphism $f$ of the scheme.

\epk

\section{Finite \'etale covers}\label{s3}
\bpk Let $\X(\F)$ be a geometrically connected variety over a number field $\kk$ which we also will treat as the sort in the model $(\X(\F), \F_{/\kk})$ for a theory $T^A$ corresponding to a geometrically connected reduced $\kk$-scheme $A.$

We develop below a uniform model for the category of all finite Galois \'etale covers of $\X.$  
\epk

\bpk \label{newnotation}
The sort $\F$ will be considered along with the projective space $\PP^N(\F),$ where $N$ is chosen so that $\X$ and \'etale covers $\Y$ of $\X$ embed in $\PP^N$ (for $\X$ a curve, $N=3$ suffices).
 
Each finite
\'etale covering of a scheme $\X$ over $\kk$ according to \ref{schemesPP}
corresponds to  a finite family of varieties 
$$\X_\mu(\F):=\{ \X_{\mu,\alpha}(\F): \mathbf{f}_\mu(\alpha)=0\}$$
and regular covering maps (morphisms) 
$$\mu:=\{ \mu_\alpha: \X_{\mu,\alpha}(\F)\to \X(\F):\ \mathbf{f}_\mu(\alpha)=0\}$$
where  $\X_{\mu,\alpha}(\F)$ are varieties defined over $\kk[\alpha],$ 
$\mathbf{f}_\mu\in \kk[T]$ is the minimal polynomial of $\alpha.$ Also 
$\mu$ is the name marking the \'etale cover
 and the collection of all such names will be denoted $\mathcal{M}_\X$
(or just   $\mathcal{M}$ if an $\X$ is clear from the context) In order to include the identity \'etale cover $1: \X\to \X$ we assume $1\in \mathcal{M}$ and  $\X_1=\X.$ 

We assume that  $\mathcal{M}_\X$ represents all finite \'etale covers of 
$\X,$ up to isomorphisms.

It is not true in general that, given $\kk(\alpha)$ and $\X_{\mu,\alpha}$  the map $\mu_\alpha$ is the only \'etale cover. Our notation just fixes one among possibly many. 

Note that $\X_{\mu,\alpha}(\F)$ in this setting is a $\kk(\alpha)$-definable subset of $\mathbf{P}^N(\F)$ and $\mu_\alpha$ a $\kk(\alpha)$ definable map, so both are given by certain formulas in the language of fields.

\medskip

Claim. We can choose family $\mu$ to be invariant under   automorphisms of $\F$ over $\kk.$ Equivalently, we may assume
$$\mu=\{ \mu_\alpha: \alpha\in \mathrm{Zeros}(\mathbf{f}_\mu)\}$$
is invariant under $\mathrm{Gal}_\kk.$ 

Proof. Given $\alpha$ and $\mu_\alpha: \X_{\mu,\alpha}(\F)\to \X(\F)$
define, for $\alpha'\in \mathrm{Zeros}(\mathbf{f}_\mu),$ $\X_{\mu,\alpha'}(\F)$ and $\mu_{\alpha'}: \X_{\mu,\alpha}(\F)\to \X(\F)$ by applying an automorphism $\sigma\in \Aut\, \F_{/\kk}$ to parameters in the definitions of $\X_{\mu,\alpha}(\F)$ and $\mu_{\alpha'}.$ Equivalently, replacing $\alpha$ by $\alpha'$ in the definitions.  


\epk

\bpk \label{numu} {\bf Intermediate morphisms between \'etale covers of $\X$.}

Given \'etale covers $\mu_\alpha: \X_{\mu,\alpha}\to \X$ and
$\nu_\beta: \X_{\nu,\beta}(\F)\to \X(\F)$ 
we consider now  \'etale morphisms
\be\label{f*} f^*:  \X_{\mu,\alpha}\to  \X_{\nu,\beta}\ee
such that $$\nu_\beta\circ f^*=\mu_\alpha.$$ 
We call such $f^*$ an intermediate morphism.
 By \ref{morphisms} and \ref{schemesPP} $f^*$ is uniquely determined by a morphism $f$ between respective schemes and the correspondence is functorial.
 
For each such $\mu_\alpha$ and $\nu_\beta$ such that a morphism (\ref{f*}) exists  we distinguish one particular morphism which we call
$$(\nu_\beta^{-1}\mu_\alpha): \X_{\mu,\alpha}\to \X_{\nu,\beta},$$
 
The notation indicates the fact that
$$\nu_\beta\circ (\nu_\beta^{-1}\mu_\alpha)=\mu_\alpha$$
not to be confused with  $\nu_\beta^{-1}\circ\mu_\alpha$ which in general is not well-defined as function but is a correspondence. 

 In general one can not make a $\mathrm{Gal}_\kk$-invariant choice of all the distingushed $(\nu_\beta^{-1}\mu_\alpha).$ 
\epk
\bpk \label{Kmuk} {\bf More notations and definitions.}
In what follows $\K_{\mu,\alpha}$ will stand for the field of definable functions of $\X_{\mu,\alpha}(\F),$ which is the same as the function field of the variety 
$\X_{\mu,\alpha}(\F)$ over $\kk[\alpha].$ Note that since the covers are Galois, $\kk[\alpha]=\kk[\alpha']$ for $\alpha'$ Galois-conjugated to $\alpha.$ Moreover, $\kk[\alpha]\cong \kk^\mu,$ the field of definable constant functions on  $\X_{\mu,\alpha}(\F).$ (Note the notation  $\kk^\mu$ which agrees with the notation $\kk^A$ of \ref{morphisms}).

It follows also that  $$\K_{\mu,\alpha}\cong \K_{\mu,\alpha'}$$
as fields. One must be aware that the isomorphism between the two fields is not canonical, however, any two isomorphisms differ by an automorphism of the field fixing the subfield $\K.$

We sometimes use the notation $$\K_\mu:=\K_{\mu,\alpha}$$ identifying the field as an abstract extension of $\K.$
\epk
\bpk {\bf Pure syntactic morphisms.}
A morphism $$(\nu_\beta^{-1}\mu_\alpha):
 \X_{\mu,\alpha}(\F)\to \X_{\nu,\beta}(\F)$$ is pure syntactic  if $(\nu_\beta^{-1}\mu_\alpha)$ as a map is a bijection of varieties and the field $\K_\mu$ of definable functions of $\X_{\mu,\alpha}(\F)$ is of the form
$\K_\nu\otimes_{\kk^\nu} \kk^\mu,$ i.e. obtained by adjoining a constant co-ordinate function with value $\alpha.$

It follows from definitions that pure syntactic $(\nu_\beta^{-1}\mu_\alpha)$ is  finite \'etale.

\medskip

Next we are going to prove that  pure syntactic covers of a given $\X_{\nu,\beta}$ can be amalgamated and the amalgamation  is pure syntactic. For this purposes we may assume consider just covers of $\X(\F).$

It is clear that up to isomorphism pure \'etale morphism $\mu_\alpha$ is determined by $\alpha\in \kk^{alg}$ so we may write it as $\alpha: \X_{\alpha}(\F)\to \X(\F).$

\medskip

{\bf Lemma.} {\em Let  $\alpha: \X_{\alpha}(\F)\to \X(\F)$ and $\beta: \X_{\beta}(\F)\to \X(\F)$ be pure syntactic. There exists a pure syntactic $\gamma: \X_{\gamma}(\F)\to \X(\F)$ amalgamating both.}

{\bf Proof.} 
Treating the schemes locally we may assume that $\X,$ $\X_\alpha$ and $\X_\beta$ are affine and their co-ordinate rings are $R,$ $A$ and $B,$ respectively. Moreover,
by definition, $\kk^R=\kk$ and $A=R\otimes_\kk\kk[a],$  $B=R\otimes_\kk \kk[b]$  
for some $a\in A$ and $b\in B$ such that
 $\kk[a]\cong \kk[\alpha]$ and  $\kk[b]\cong \kk[\beta].$ 
 
 Moreover,
 $a(x)=\alpha$ on $\X_\alpha(\F)$ and $b(x)=\beta$ on $\X_\beta(\F).$

Let $\kk[\gamma]=\kk[\alpha,\beta].$ Let the ring
$C\cong_\phi R\otimes_\kk \kk[\gamma],$ where $\gamma$ corresponds to an abstract element $c\in C$ by the isomorphism $\phi.$ Let $a'$ and $b'$ be the uniquely defined elements of $C$ which correspond under $\phi$ to $\alpha$ and $\beta$ respectively. 

We now have obvious embeddings 
$$g: A\to C,\ a\mapsto a'\mbox{ and } 
h: B\to C,\ b\mapsto b'.$$

These give us \'etale morphisms of respective schemes and we can now construct a model of $T^C$ given by a variety
$\X_\gamma(\F)$ which is set-wise just $\X(\F)$ and which has constant co-ordinate functions 
$a'(x)=\alpha$ and $b'(x)=\beta$ for $x\in \X(\F).$ By our choices $\X_\alpha(\F)$ and $\X_\gamma(\F)$ is an appropriate pair of models of 
$T^A$ and $T^C$ and so is the pair $\X_\beta(\F)$ and $\X_\gamma(\F).$ Thus we have defined morphisms
$$g^*_{\gamma,\alpha}: \X_\gamma(\F)\to \X_{\alpha}(\F)\mbox{ and } h^*_{\gamma,\beta}: \X_\gamma(\F)\to \X_{\beta}(\F),$$ as required. $\Box$ 

\epk

\bpk \label{pureG} {\bf Pure geometric  morphisms.}
A morphism $$(\nu_\beta^{-1}\mu_\alpha): \X_{\mu,\alpha}(\F)\to \X_{\nu,\beta}(\F),$$ is pure geometric (semantic)  if the field $\kk^\mu$ of definable constant functions of $\X_{\mu,\alpha}(\F)$ is equal to $\kk^\nu.$ That is both varieties are defined over $\kk[\beta]$ and $\kk[\beta]=\kk[\alpha].$

{\bf Warning.} It is not true in general that 
$(\nu_\beta^{-1}\mu_\alpha)$ is defined over $\kk[\beta].$

\medskip

{\bf Remark.} If $(\nu_\beta^{-1}\mu_\alpha)$ is pure geometric then
both absolutely irreducible varieties $\X_{\mu,\alpha}(\F)$ and $\X_{\nu,\beta}(\F)$ are defined over $\kk[\beta]$ and hence their function fields, as varieties over $\F,$ are $\K_{\mu}\otimes_{\kk[\beta]} \F$ and $\K_{\nu}\otimes_{\kk[\beta]} \F$ respectively.
 
Since $\kk[\beta]$ is algebraically closed in each of the fields,  
\be\label{Deck0} \mathrm{Gal}(\K_{\mu}\otimes_{\kk[\beta]} \F : \K_\nu\otimes_{\kk[\beta]} \F)\cong\mathrm{Gal}(\K_\mu:\K_\nu). \ee
\epk
\bpk \label{geomcover} {\bf Lemma.} {\em Every finite   \'etale $\mu_\alpha: \X_{\mu,\alpha}(\F)\to \X(\F)$ can be presented as a composition of pure syntactic
 $\nu_\beta: \X_{\nu,\beta}(\F)\to \X(\F)$ and a pure geometric $(\nu_\beta^{-1}\mu_\alpha): \X_{\mu,\alpha}(\F)\to \X_{\nu,\beta}(\F).$  
}

{\bf Proof.} First assume $\X$ to be affine and let $A$ to be its coordinate ring over $\kk.$ Define the coordinate ring of $\X_{\nu,\beta}$ to be 
$\kk^\mu\otimes_\kk A.$ 
Now the embedding $A\hookrightarrow \kk^\mu\otimes_\kk A$ induces a pure syntactic morphism   $\nu_\beta: \X_{\nu,\beta}(\F)\to \X(\F)$ which further lifts to pure geometric
 $(\nu_\beta^{-1}\mu_\alpha): \X_{\mu,\alpha}(\F)\to \X_{\nu,\beta}(\F).$  
 
 The same procedure applies to open charts of $\X$ in the general case. $\Box$ 
\epk
\bpk In the following chapters we will be interested in the special cases of complex algebraic varieties. This corresponds to the case $\F=\C.$  In particular, if  $(\nu_\beta^{-1}\mu_\alpha)$
is pure geometric then by Remark in \ref{pureG} $(\nu_\beta^{-1}\mu_\alpha)$ is a morphism in the category of complex algebraic varieties, and hence in the category of complex manifolds.

We are going to use a standard fact in the theory of \'etale covers (easily following from definitions):  {\em a morpism between complex algebraic varieties is \'etale if and only if it is unramified in the category of complex manifolds}.

\epk

\bpk \label{defDeck}
Recall the definition of the deck-transformation group  for an unramified cover:

$\mathrm{GDeck}(\X_{\mu,\alpha}(\C)/\X_{\nu,\beta}(\C))$ is the group of biholomorphic transformations of $\X_{\mu,\alpha}(\C)$ as a complex manifold
 fixing fibres  of $\X_{\mu,\alpha}(\C)\to \X_{\nu,\beta}(\C).$ 
 
 By standard algebra/analytic comparison facts we get the same group by replacing  ``biholomorphic transformations of $\X_{\mu,\alpha}(\C)$ as a complex manifold'' by  ``biregular transformations of $\X_{\mu,\alpha}(\C)$ as a complex algebraic variety.

 This can be equivalently given in terms of function fields of the varieties defined over respective subfieds, 
$$ \mathrm{GDeck}(\X_{\mu,\alpha}(\C)/\X_{\nu,\beta}(\C))\cong
\mathrm{Gal}(\K_{\mu}\otimes_{\kk[\alpha]} \C : \K_\nu\otimes_{\kk[\beta]} \C).$$ 
In the right-hand side  we can replace $\C$ by any algebraically closed field $\F$ containing $\kk[\alpha]$ and $\kk[\beta].$

Given  pure geometric morphism $$(\nu_\beta^{-1}\mu_\alpha): \X_{\mu,\alpha}(\F)\to \X_{\nu,\beta}(\F),$$
we call the {\bf geometric deck-transformation group} the group of biregular transformations of $\X_{\mu,\alpha}(\F)$  fixing fibres  of $\X_{\mu,\alpha}(\F)\to \X_{\nu,\beta}(\F)$ as  algebraic varieties over $\F.$

We often refer to a  deck-transformation group without mentioning $\F$ and we always assume
 that the cover is pure geometric. In particular, one easily  deduces from the above, 

\be\label{GDeck} \mathrm{GDeck}(\X_{\mu,\alpha}/\X_{\nu,\beta})\cong \mathrm{Gal}(\K_\mu:\K_\nu). \ee

\epk

\bpk \label{rem2} {\bf Lemma.} {\em  Suppose   $(\nu_\beta^{-1}\mu_\alpha)': \X_{\mu,\alpha}(\F)\to \X_{\nu,\beta}(\F)$
is a pure geometric morphism with the property $\nu_\beta\circ(\nu_\beta^{-1}\mu_\alpha)'=\mu_\alpha.$ 
 Then
there is $g\in \mathrm{GDeck}(\X_{\nu,\beta}/\X)$ such that $$(\nu_\beta^{-1}\mu_\alpha)'=g\circ (\nu_\beta^{-1}\mu_\alpha).$$}

{\bf Proof.} Consider $$S:=\{ \la x,x'\ra\in \X_{\nu,\beta}(\F)\times  \X_{\nu,\beta}(\F): \exists y\in   \X_{\mu,\alpha}\,
(\nu_\beta^{-1}\mu_\alpha)(y)=x\ \&\ (\nu_\beta^{-1}\mu_\alpha)'(y)=x'\}.$$
We may assume $\F=\C$ and work in the category of complex manifolds.

Since \'etale coverings are closed unramified maps, $S$ is a Zariski closed subset of the set
$$E:=\{ \la x,x'\ra\in \X_{\nu,\beta}(\F)\times  \X_{\nu,\beta}(\F): \nu_\beta(x)=\nu_\beta(x')\}.$$
and is locally holomorphically isomorphic to  $\X_{\nu,\beta}(\F).$

We claim that the irreducible components of $E$ are in one-to-one correspondence with elements
$g\in \mathrm{GDeck}(\X_{\nu,\beta}/\X)$ and are indeed
of the form
$$E_g=\{   \la x,x'\ra: g*x=x'\}.$$
Indeed, pick a point $\la x_0,x'_0\ra\in E$ and an irreducible   component $C_0$ of $E$ containing
the point. Since $\mathrm{GDeck}(\X_{\nu,\beta}/\X)$ acts transitively on fibres of $\nu_\beta$ there is $g$ such that $g*x_0=x'_0,$ that is  $\la x_0,x'_0\ra\in E_g.$ Thus $C_0$ and $E_g$ have a common point which is only possible if $C_0=E_g,$ since $E$ is non-singular.

The claim implies the required. 

\epk
\bpk \label{projlimit} {\bf Projective limit of \'etale covers along distinguished morphisms.}
The scheme-theoretic category of finite Galois \'etale covers allows a profinite universal object, the fundamental pro-object in terminology of SGA1, see also \cite{MilneEC}, pp.26-27.  Our construction below in this section corresponds to a different object which depends on  extra data, namely a system of distinguished morphisms. The structure which properly represents Grothendieck's pro-object will be constructed later in section  \ref{s5+}.

We work with the family $$\mathcal{X}:= \{  \X_{\mu,\alpha}(\F): \ \mu\in \MM_X, \alpha\in \mathrm{Zeros}\,\mathbf{f}_\mu\}$$
of varieties, covers of $\X,$ and the family of all the intermediate morphisms between them.
   
We 
aim to construct a projective system of 
distinguished intermediate morphisms $(\nu_\beta^{-1}\mu_\alpha): \X_{\mu,\alpha}(\kk^{alg})\to \X_{\nu,\beta}(\kk^{alg}),$ equivalently, distinguished morphisms $(\nu_\beta^{-1}\mu_\alpha): \X_{\mu,\alpha}(\F)\to \X_{\nu,\beta}(\F),$
satisfying the conditions:

(i) For each $\mu_\alpha$
a distinguished $(\mu_\alpha^{-1}\mu_\alpha): \X_{\mu,\alpha}\to  \X_{\mu,\alpha}$ is an identity on  $\X_{\mu,\alpha}.$

(ii) if $(\nu_\beta^{-1}\mu_\alpha): \X_{\mu,\alpha}(\F)\to \X_{\nu,\beta}(\F)$ and $(\lambda_\gamma^{-1}\nu_\beta): \X_{\lambda,\gamma}(\F)\to \X_{\mu,\alpha}(\F)$ are distinguished then their composite $$(\lambda_\gamma^{-1}\mu_\alpha):=(\nu_\beta^{-1}\mu_\alpha)\circ (\lambda_\gamma^{-1}\nu_\beta)$$
is distinguished.

In the Appendix (section \ref{Ap}) we prove:

{\bf Claim}.{\em  A projective system of distiguished morphisms exists.}

Hence there is a well-defined projective limit  
$$\tilde{\X}^{et}(\F):=\lim_{\leftarrow} \X_{\mu,\alpha}(\F)$$
along the system of morphisms $(\nu_\beta^{-1}\mu_\alpha)$
 
This gives us a cover of $\X(\F),$
$${\pp}:\tilde{\X}^{et}(\F)\to \X(\F)$$
and of each intermediate $\X_{\mu,\alpha}(\F),$  
$$\pp_{\mu,\alpha}:\tilde{\X}^{et}(\F)\to \X_{\mu,\alpha}(\F)$$
which satisfy the relations
\be \label{(d)}  (\nu_\beta^{-1}\mu_\alpha)\circ\pp_{\mu,\alpha}=\pp_{\nu,\beta}\ee
for the distingushed $(\nu_\beta^{-1}\mu_\alpha).$




\epk

\bpk \label{strings} {\bf Remarks.} 

1. The construction of $\X^{et}$ depends on the choice of the projective system of distinguished morphisms.

2. By definition
$$\tilde{\X}^{et}(\F)\subset (\PP^N(\F))^{\MM_\X}$$
where the right-hand side is the set of all functions  $u:\MM_\X\to \PP^N(\F).$ 

The automorphism group $\Aut\, \F_{/\kk}$ acts on 
$\PP^N(\F).$ In general, $\tilde{\X}^{et}(\F)$ is not  $\Aut\, \F_{/\kk}$-invariant. 

In particular, it is not true in general that the families $$\pp_\mu:=\{ \pp_{\mu,\alpha}: \alpha\in \mathrm{Zeros}\, \mathbf{f}_\mu\}$$ are
 $\Aut\, \F_{/\kk}$ -invariant  even if $\F=\kk^{alg}.$ However, the families of sets
$$\tilde{\pp}_{\mu}\ := \{ g\circ\pp_{\mu,\alpha}: \ g\in \mathrm{GDeck}(\X_{\mu,\alpha}/\X), \alpha\in \mathrm{Zeros}\, \mathbf{f}_\mu\}$$
are $\Aut\, \F_{/\kk}$-invariant since the initial system $\{ \mu_\alpha:   \alpha\in \mathrm{Zeros}\, \mathbf{f}_\mu\}$ is  $\Aut\, \F_{/\kk}$-invariant.

\epk

\bpk \label{warning} {\bf Remark.} 
The  maps
$\pp_{\mu,\alpha}$ are determined by the $\mu_\alpha$ up two regular transformation of $\X_{\mu,\alpha}$ preserving fibres of $\mu_\alpha,$ that is up to the action of $g\in \mathrm{GDeck}(\X_{\mu,\alpha}/\X).$ Indeed,
 $$\mu_\alpha=  \pp\circ (g\circ\pp_{\mu,\alpha})\inv.$$

To see this we note that $\pp_{\mu,\alpha}$ is the inverse limit of the system of distinguished $(\mu_\alpha^{-1}\nu_\beta),$ for $\nu\in \MM_\X,$ ``above'' $\mu,$ so is 
 uniquely determined by the choices of $(\mu_\alpha^{-1}\nu_\beta): \X_{\nu,\beta}\to \X_{\mu,\alpha}.$ A different choice $\pp'_{\mu,\alpha}$ results from a different choice
  $(\mu_\alpha^{-1}\nu_\beta)': \X_{\nu,\beta}\to \X_{\mu,\alpha}$ at some stage  $\nu\in \MM_\X.$ By \ref{rem2}
$$(\mu_\alpha^{-1}\nu_\beta)'=g\circ (\mu_\alpha^{-1}\nu_\beta),\mbox{ for some }g\in \mathrm{GDeck}(\X_{\mu,\alpha}/\X)$$
which proves the claim.

\epk
\bpk \label{pmu} {\bf Lemma.} {\em Suppose $\X_{\mu,\alpha}(\kk^{alg})=\X_{\nu,\beta}(\kk^{alg})$ as sets and the set-to-set maps
$$\mu_\alpha:  \X_{\mu,\alpha}(\kk^{alg})\to \X(\kk^{alg})\mbox{ and }\nu_\beta:\X_{\nu,\beta}(\kk^{alg}) \to \X(\kk^{alg})$$
are equal. Then $$\pp_{\mu,\alpha}=\pp_{\nu,\beta}.$$
}

{\bf Proof.} The assumptions of the lemma imply that $\X_{\mu,\alpha}(\kk^{alg})=\X_{\nu,\beta}(\kk^{alg})$ as varieties over $\kk^{alg}$ and $\mu_\alpha=\nu_\beta$ as morphisms over $\kk^{alg}.$ Hence $\pp_{\mu,\alpha}=\pp_{\nu,\beta}$  by construction. $\Box$

\epk

\bpk
Along with the projective limit $\tilde{\X}^{et}(\F)$ of covers one obtains also the projective limit of deck transformation
groups which we denote

$$\hat{\Gamma}:=\lim_{\leftarrow} \mathrm{GDeck}(\X_{\mu,\alpha}/\X).$$

This group acts freely on $\tilde{\X}^{et}(\F)$ in agreement with the actions of  $\mathrm{GDeck}(\X_{\mu,\alpha}/\X)$
on $\X_{\mu,\alpha}(\F).$

In terms of the action of $\hat{\Gamma}$ on $\tilde{\X}^{et}(\F)$ one defines period subgroups 
$$\hat{\Delta}_{\mu,\alpha}:=\mathrm{Per}(\pp_{\mu,\alpha}):=
\{ \gamma\in \tilde{\Gamma}: \forall x\ \pp_{\mu,\alpha}(\gamma*x)=\pp_{\mu,\alpha}(x)\}.$$

It follows from \ref{geomcover} and (\ref{GDeck})
 that  $\hat{\Delta}_\alpha$ is a normal subgroup of $\hat{\Gamma}$ and that
\begin{itemize}
\item $\mathrm{GDeck}(\X_{\mu,\alpha}/\X)\cong \hat{\Gamma}/\hat{\Delta}_{\mu,\alpha}.$
\end{itemize}

In particular, $\hat{\Delta}_{\mu,\alpha}$ is of finite index in $\hat{\Gamma}.$ Moreover, one sees that $$\bigcap_{{\mu,\alpha}}\hat{\Delta}_{\mu,\alpha}=\{ 1\},$$ that is 
$\hat{\Gamma}$ is residually finite.

\epk

\bpk {\bf The language $\mathcal{L}_\X^\sharp.$}

We consider $\tilde{\X}^{et}(\F)$ in a formalism of a structure with two basic sorts  $\tilde{\X}^{et}(\F)$ and $ \F_{/\kk}.$ The sorts interact via the families $\pp_{\mu,\alpha}$ of maps which we formalise as follows.

We use $\pp_\mu,$ $\mu\in \MM_\X,$ as names which are interpreted as  maps 
$$\pp_\mu:\tilde{\X}^{et}(\F)\times \mathrm{Zeros}\,\mathbf{f}_\mu \to \PP^N(\F)$$
such that $$\pp_\mu(u,\alpha):=\pp_{\mu,\alpha}(u),$$
where on the right we use the notation of previous subsections for individual $\alpha$-definable maps. (So we continue to think of $\pp_\mu$ as families of maps.)

We will denote 
$$\XX^{et}_\sharp(\F):= \left( \tilde{\X}^{et}(\F),\, \F_{/\kk}, \{ \pp_\mu\}_{\mu\in \mathcal{M}_\X}\right)$$ 
the structure obtained by the inverse limit construction from the $\F$-model $\X(\F)$ of the curve $\X,$  where $\F$ is an algebraically closed field containing $\kk.$

The language of the structure: sort $\U$ for $\tilde{\X},$ sort $\F$ with the language of rings, and names (for the families) $\pp_\mu$ of maps, we denote as
 $\mathcal{L}_\X^\sharp.$

The superscript $\sharp$ indicates that the language  we use here for describing the universal \'etale cover is {\em excessive}, and thus not all the possible symmetries of the cover can be realised as automorphisms. The adequate language $\mathcal{L}_\X$ and the adequate \'etale cover structure ${\XX^{et}}$ will be introduced in section  \ref{s5+}.

\medskip

We now note that the language $\mathcal{L}_\X^\sharp$ of the  structure is sufficient to express all the notions of \ref{newnotation} - \ref{projlimit}. Namely,
$$\X_{\mu,\alpha}(\F)=\pp_{\mu,\alpha}(\tilde{\X}(\F))$$
as a subset of $\PP^N(\F).$ Then the $\kk(\alpha)$-definable structure on $\X_{\mu,\alpha}(\F)$ is defined by the embedding in $\PP^N$ over $\F_{/\kk}.$

\epk

\bpk\label{tilde0} {\bf Remark.}
 Each $v\in \tilde{\X}(\F)$  can be identified with the
type $\tau_v$ in one variable $u$ 
$$\tau_v=\{ u\in \U\ \& \ \pp_{\mu,\alpha}(u)=\pp_{\mu,\alpha}(v): \ \mu\in \mathcal{Mu},\, \mathbf{f}_\mu(a)=0   \}.$$ 
This is a type over the  countable subset of $\F$

 $$\mathrm{Supp}(v)=\{ \pp_{\mu,\alpha}(v): \ \mu\in \mathcal{Mu},\, \mathbf{f}_\mu(a)=0   \}.$$
\epk

\section{Analytic Covers }\label{sAP}
In this section we are using the language $\mathcal{L}_\X^\sharp$ to describe the analytic cover $\tilde{\X}^{an}(\C)$ of the complex variety 
$\X(\C),$ seen also as a complex manifold.
\bpk \label{complexX} We start agin with a reduced smooth geometrically connected quasi-projective  $\kk$-scheme $\X$ and $\X(\C)\subset \PP^N(\C)$ as in \ref{defV} and section \ref{s3}.


\epk
\bpk \label{GP} Let $$\pp: \tilde{\X}^{an}(\C)\to \X(\C)$$ be the  universal cover (complex manifold) of the complex variety $\X(\C)$ 
with the covering holomorphic map. The (topological) fundamental group $\Gamma$ of $\X(\C)$ acts properly discontinuously  on $\tilde{\X}^{an}(\C),$ so
\be \label{gammaX}\Gamma \backslash \tilde{\X}^{an}(\C)\cong \X(\C)\ee
as complex analytic manifolds.


\medskip

It is clear that the cartesian product $\Gamma^n$ acts properly
discontinuously on the $n$-th cartesian power of 
$\tilde{\X}^{an}(\C).$
\epk

\bpk Let $\mu_\alpha: \X_{\mu,\alpha}(\C)\to \X(\C)$ be a finite  regular covering map as in section \ref{s3} which corresponds to an \'etale covering. Then $\mu_\alpha$ is unramified in the sense of complex Hausdorff topology.
 Since $\tilde{\X}^{an}(\C)$  is the universal cover of $\X(\C)$
  there is a normal subgroup
$\Delta_{\mu,\alpha}$ of $\Gamma$ of finite index and a holomorphic map
 $$\pp_{\mu,\alpha}: \tilde{\X}^{an}(\C)\to \X_{\mu,\alpha}(\C)\cong {\Delta}_{\mu,\alpha}\backslash\ \tilde{\X}^{an}(\C)$$
 where the last isomorphism is understood as a biholomorphic isomorphism between complex manifolds. 
 It is clear from the general facts that the group $\Gamma/\Delta_{\mu,\alpha}$ acts on $\X_{\mu,\alpha}(\C)$ and 
 $$ \Gamma/\Delta_{\mu,\alpha}\cong \mathrm{GDeck}(\X_{\mu,\alpha}(\C)/\X(\C)).$$
  
Using that action of this group one can always adjust the choice of $\pp_{\mu,\alpha}$ so that
for a given finite  unramified cover $\mu_\alpha: \X_{\mu,\alpha}(\C)\to \X(\C),$
 
\begin{itemize}
\item
 $\mu_\alpha\circ \pp_{\mu,\alpha}=\pp$
\end{itemize}
Note also that by definition
\begin{itemize}
\item
$\pp_{\mu,\alpha}(x)=\pp_{\mu,\alpha}(x')\Leftrightarrow \Delta_{\mu,\alpha}\cdot x=\Delta_{\mu,\alpha}\cdot x'
$
\end{itemize}

\epk
\bpk\label{an/alg} {\bf Finite unramified covers and \'etale covers of $\X(\C).$} The biholomorphic isomorphism type of  ${\Delta}_{\mu,\alpha}\backslash\ \tilde{\X}^{an}(\C)$  according to algebraic/analytic comparison theorems corresponds to the isomorphism type of the complex algebraic variety  $\X_{\mu,\alpha}(\C).$  So, if 
$\X_{\nu,\beta}(\C)$ is a pure syntactic cover of $\X_{\mu,\alpha}(\C),$ then the respective complex manifolds are the same and $\Delta_{\mu,\alpha}=\Delta_{\nu,\beta}.$

Conversely, when a normal subgroup $\Delta$ of $\Gamma$ of finite index is given one can always identify the complex manifold ${\Delta}\backslash\ \tilde{\X}^{an}(\C)$ as an unramified cover of the complex manifold $\X(\C)$ and by the Riemann Existence Theorem  ${\Delta}\backslash\ \tilde{\X}^{an}(\C)$ can be identified as an algebraic variety over $\C,$ \'etale covering $\X(\C).$ Since $\X$ is defined over
 $\kk,$ the algebraic variety ${\Delta}\backslash\ \tilde{\X}^{an}(\C)$ can be defined over an algebraic extension
 $\kk(\alpha)$ of $\kk,$ thus taking the form  $\X_{\mu,\alpha}(\C).$

\epk
\bpk {\bf The standard analytic structure.} 
The two-sorted structure  
$$\XX^{an}_\sharp(\C)= \left( \tilde{\X}^{an}(\C),\,  \C_{/\kk}, \{ \pp_\mu\}_{ \mu\in \mathcal{M}_\X} \right) $$
where, as above, $\tilde{\X}^{an}(\C)$ is the complex universal cover of $\X(\C)$ seen as a set, 
$ \C_{/\kk}$ is the complex numbers in the language of rings and names for points of $\kk.$ 
For each $\mu\in \mathcal{M}_\X$ and each zero $\alpha$ of corresponding polynomial $\mathbf{f}_\mu$ there is a $\pp_{\mu,\alpha}\in \pp_\mu,$
$$\pp_{\mu,\alpha}:\tilde{\X}^{an}(\C) \to \X_{\mu,\alpha}(\C)\subset \PP^N(\C).$$

\epk
\section{The first order theory of  covers} \label{s6}
\bpk \label{Tsharp} {\bf The theory $T^\sharp_\X$}

The axioms describe a two-sorted structure $(\U;\F, \{ \pp_\mu\}_{\mu\in \mathcal{M}_\X})$:
\begin{itemize}
\item[A1] {\em $\F$ is an algebraically closed field of characteristic $0$ with subfield $\kk$ of constants.}

Remark. $\X(\F)$ and each of the varieties 
$\X_{\mu,\alpha}(\F)$ together with \'etale morphisms $\mu_\alpha$ and $(\nu_\beta^{-1}\mu_\alpha)$ are given by $\kk^{alg}$-definable subsets of $\PP^N(\F)$ and $\kk^{alg}$-definable maps between them, as described in \ref{newnotation} - \ref{projlimit}.

\item[A2]($\mu$) {\em Given $\mu\in \mathcal{M}_\X,$ 
 $\pp_{\mu}$ is a map  with domain $\U\times \mathrm{Zeros}\,\mathbf{f}_\mu$ and values in $ \PP^N(\F),$ and
 $$\pp_\mu(\U\times \{ \alpha\})= \X_{\mu,\alpha}(\F).$$}

\item[A3]($\mu,\nu$) {\em There is an evaluation of parameters $\alpha,\beta,\ldots $ in $\kk^{alg}\subseteq \F$ such that for all  
 $\mu,\nu\in \mathcal{M},$ 
 for respective zero $\alpha$ of $\mathbf{f}_\mu$ and $\beta$ of $\mathbf{f}_{\nu}$}

$$(\nu_\beta^{-1}\mu_\alpha)\circ \pp_{\mu,\alpha}=  \pp_{\nu,\beta}.$$ 

Remark. A3 is given by an infinite system of $\exists \forall$-sentences each of which states that a finite system of identities (\ref{(d)}) which hold in $\XX^{et}_\sharp(\F)$ is consistent.


\end{itemize}

\epk
\bpk \label{Comparison}
{\bf Comparison Theorem.} {\em $\XX^{an}_\sharp(\C)$ and $\XX^{et}_\sharp(\F)$
 are models of $T^\sharp_\X,$ with $\U=\tilde{\X}^{an}(\C)$ and $\F=\C$ in the first case and $\U=\tilde{\X}^{et}(\F)$ in the second case.}

 {\bf Proof.} Immediate from \ref{newnotation} - \ref{projlimit}. $\Box$
 
\epk

\bpk \label{sharpk} Let $G^\sharp$ be the subgroup of  $\mathrm{Gal}_\kk$ which preserves  the family of all the distinguished morphisms $(\nu_\beta^{-1}\mu_\alpha),$ that is,

  $\sigma\in G^\sharp$  
 if and only if 
 $(\nu_\beta^{-1}\mu_\alpha)^\sigma=({\nu_{\sigma(\beta)}^{-1}}\mu_{\sigma(\alpha)}),$ distinguished.

\medskip

Define $\kk^\sharp$ to be the subfield of $\kk^{alg}$ which is point-wise fixed by $G^\sharp.$

We note that in $\XX^{et}(\F)$ by construction 
 $\sigma(\pp_{\mu,\alpha})=\pp_{\mu,\sigma(\alpha)}$ for $\sigma\in G^\sharp,$ since
  the $\pp_{\mu,\alpha}$ are limits of chains of distinguished morphisms.
  
 Consider the sort $\F$ with 
the families of the distinguished morphisms as a structure, call it $\F^\sharp_{\kk}.$ 
\epk
\bpk\label{sharpLem} {\bf Lemma.} {\em In $\F^\sharp_{\kk}$ 
$$\kk^\sharp= \F\cap \dcl (\emptyset).$$ 

In any model of $T^\sharp_\X$
  $$\kk^\sharp\subseteq \F\cap \dcl (\emptyset).$$}

{\bf Proof.}  $\F^\sharp_{\kk}$  is interpretable in 
 the  field structure $\F$ with parameters. By elimination of imaginaries in the theory of algebraically closed fields,  $\F^\sharp_{\kk}$ is bi-interpretable with $\F_{\kk^\sharp}.$ 
 
The second statement is just a consequence of the first one. 
 $\Box$

\epk
The next theorem was proved by A.Harris \cite{Adam} in a less flexible language.

\bpk\label{complete} {\bf Theorem.}

1. { \em The first-order theory $T^\sharp_\X$ is complete. In particular,
$$\XX^{an}_\sharp(\C)\equiv \XX^{et}_\sharp(\F)$$
}

2. {\em In a model $(\U, \F_{/\kk})$ of the theory,
a $\mathcal{L}_\X^\sharp$-first-order definable subset of $\U^n$ is of a form
$\pp\inv(S)$
for  $S\subseteq \X^n(\F),$ $\mathcal{L}_\X^\sharp$-definable.

A first order $\mathcal{L}_\X^\sharp$-definable subset of $\X^n(\F)$ is definable in the field language using parameters in $\kk^\sharp.$ $T^\sharp_\X$ has elimination of quantifiers in the language $\mathcal{L}^\sharp_\X$ expanded by names for $\kk^\sharp.$}

\medskip

3. {\em $T^\sharp_\X$ is superstable. }

 {\bf Proof.} We may assume by \ref{sharpLem} that $\mathcal{L}^\sharp_\X$ has names for elements of $\kk^\sharp,$ so below $\F$ is always
 $\F_{/\kk^\sharp}.$
  
  Let $(\U, \F)$ be a $\kappa$-saturated model $T^\sharp_\X,$ $\kappa$ uncountable cardinal. For each $u\in \U$ define the type $\tau_u$ in variable $v$ over $\F$
 $$\tau_u(v):=\{ \pp_{\mu,\alpha}(v)=\pp_{\mu,\alpha}(u)\}$$

We will construct an elementary submodel
 $$(\U^*, \F^*)\preccurlyeq (\U, \F)$$
 such that
 \begin{itemize}
 \item[(i)] $|\U^*|=| \F^*|=\kappa.$
 \item[(ii)] any type of the form $\tau_u,$ for $u\in \tilde{\X}( \F^*),$ is realised in $\U^*$ by exactly $\kappa$ distinct elements.
 \item[(iii)] any element of  $\U^*$ realises a type  $\tau_u,$ for $u\in \tilde{\X}( \F^*).$
\end{itemize}  
We call $(\U^*, \F^*)$ as above {\bf a $\kappa$-good elementary submodel of } $(\U,\F).$

(Remark. Any $\kappa$-saturated model of cardinality $\kappa\ge \mathfrak{c}$ is $\kappa$-good.
 Saturated model of  cardinality $\kappa\ge \mathfrak{c}$ exist provided CH holds or $T_\X$ is stable.)

\medskip

Let $\F^0\subseteq \F$ be an algebraically closed subfield of cardinality $\kappa.$ By axiom A2 and A3 each type $\tau_u$ is realised in a saturated enough model of $T^\sharp_\X,$ so
we can embed $$\tilde{\X}(\F^0)\subseteq \U.$$
For each $u\in  \tilde{\X}(\F^0)$ the set  $\tau_u(\U)$ of realisations of the type $\tau_u$  in $\U$ is at least of cardinality $\kappa.$  Let $S^0_u\subseteq \tau_u(\U)$ be a subset of cardinality exactly $\kappa.$ 

Let $$\U^0=\bigcup \{ S^0_u: u\in  \tilde{\X}(\F^0)\}.$$
$(\U^0,\F^0)$ is a submodel of $(\U,\F)$ (check the axioms of $T_\X$) satifying (i) and (ii),
but we can not 
claim it is an elementary submodel. By L\o wenheim-Skolem we can construct 
 $$(\U^0,\F^0)\subseteq  (\U^{(0)}, \F^{(0)})\preccurlyeq (\U, \F)$$
 such that $|\U^{(0)}|=|\F^{(0)}|=\kappa.$
 
 Now we repeat our construction starting with $\F^1=\F^{(0)}$ in place of $\F^0$ and set
 $$\U^1=\bigcup \{ S^1_u: u\in  \tilde{\X}(\F^1)\},$$
 where $S^1_u\supseteq S^0_u$ for $u\in  \tilde{\X}(\F^0).$
 Again, $(\U^1,\F^1)$ is a submodel of $(\U,\F)$  satifying (i) and (ii) and we can continue
 $$(\U^0,\F^0)\subseteq  (\U^{(0)}, \F^{(0)})\subseteq  (\U^1, \F^1)\subseteq  (\U^{(1)}, \F^{(1)})\ldots \subseteq  (\U^i, \F^i)\subseteq  (\U^{(i)}, \F^{(i)})\ldots $$
 where all models satisfy (i) and (ii) and 
   $$ (\U^{(0)}, \F^{(0)})\preccurlyeq  (\U^{(1)}, \F^{(1)})\ldots   \preccurlyeq    (\U, \F)$$ 
 Set   $$(\U^*, \F^*)=\bigcup_{i<\omega}  (\U^{(i)}, \F^{(i)}).$$
 This satisfies all the requirements.

\medskip

{\bf Claim.} For any $\kappa$-good models  $(\U_1,\F_1)$ and $(\U_2,\F_2)$  of cardinality $\kappa$
  there exists an isomorphism $$\mathbf{i}:(\U_1,\F_1)\to (\U_2,\F_2)$$
  
  Proof. The fields in both structures have to be of the same cardinality $\kappa$ and hence they are isomorphic over $\kk^\sharp,$ the subfield of definable elements. We assume without loss of generality that $\F_1=\F_2=:\F,$ and $\mathbf{i}$ is an identity map on $\F.$

Now $\X$ and all the  $\X_{\mu,\alpha}$ along with morphisms
$\mu_\alpha$ have the same meaning $\X(\F),$ $\X_{\mu,\alpha}(\F)$ and so on,
in the two structures.

We need to construct $\mathbf{i}: \U_1\to \U_2.$

By assumptions each type $\tau_u$
 is realised in both models $(\U_1,\F)$ and $(\U_2,\F)$ by $\kappa$-many points of the sort $\U,$ call the set of realisations
 $\tau_u^{\U_1}$ and  $\tau_u^{\U_2},$ respectively. Moreover, 
$$\U_1=\bigcup_u  \tau_u^{\U_1}\mbox{ and }\U_2=\bigcup_u  \tau_u^{\U_2}.$$
 Set $\mathbf{i}: \U_1\to \U_2$ to be equal to $\bigcup_u \tau_u,$ where 
 $\mathbf{i}_u: \tau_u^{\U_1}\to \tau_u^{\U_2}$ are   bijections.
This preserves all the maps $\pp_{\mu_\alpha}$ and hence is an isomorphism.
Claim proved.

It follows that any two $\kappa$-saturated models of $T^\sharp_\X$ have isomorphic elementary submodels, that is the models are elementarily equivalent.
The first statement of the Theorem follows. 

\medskip

2. Follows from the claim. 

\medskip

3. Is a direct consequence of 2.
$\Box$ 
\epk

\bpk \label{presection} {\bf Proposition.} {\em 
 $$\Aut\, \XX^{et}_\sharp(\F)\cong\Aut\, \F_{/\kk^\sharp}$$
canonically. In particular,  
$$\Aut\, \XX^{et}_\sharp(\kk^{alg})\cong \mathrm{Gal}_{\kk^\sharp}$$

}

{\bf Proof.} It follows from \ref{complete}.3 and the fact that the theory of algebraically closed fields eliminates imaginaries
 that the substructure 
$\F_{/\kk^\sharp}$ on sort $\F$
 is homogeneous.  An automorphism $\sigma\in  \Aut\, \F_{/\kk^\sharp}$ induces a unique bijection on the space of complete types of the form $\tau_u$ 
  which induces a bijection on $\tilde{\X}(\F),$
   an automorphism of the structure, by  \ref{complete}.2.
  $\Box$
\epk

\section{The language of universal covers and the \'etale fundamental group}\label{s5+}

\bpk \label{ppp} Define, for each $\mu\in  \mathcal{M}_\X$ and $\alpha\in \mathrm{Zeros}\,\mathbf{f}_\mu$ the set of maps
$$\tilde{\pp}_{\mu,\alpha}:= \{ p: \U\to \X_{\mu,\alpha}|\  \ p=g\circ \pp_{\mu,\alpha},\  g\in \mathrm{GDeck}(\X_{\mu,\alpha}/\X)\}.$$

This is a $\mathrm{GDeck}(\X_{\mu,\alpha}/\X)$-set, that is there is a canonical free and transitive action of the group on the set.

Define, for each $\mu\in \MM_\X,$
 $$\tilde{\pp}_{\mu}=\bigcup\{ \tilde{\pp}_{\mu,\alpha}: \mathbf{f}_{\mu}(\alpha)=0\},$$
 the  finite set which is split into subsets  indexed by $\alpha\in \mathrm{Zeros}\,\mathbf{f}_{\mu}.$ 
\epk

 \bpk\label{Lflat}
We define the language $\mathcal{L}_\X$ to contain:
\begin{itemize}
\item sorts $\U,$ $\F$ and 
finite sorts $\tilde{\pp}_{\mu}.$ 
for
$\mu\in \mathcal{M}_\X.$ 

The subsets  $\tilde{\pp}_{\mu,\alpha}$ of sorts  $\tilde{\pp}_{\mu}$ are definable using binary predicates which
we write as a relation $p\in  \tilde{\pp}_{\mu,\alpha}$ between $p\in  \tilde{\pp}_{\mu}$ and $\alpha\in \F.$
 
\item $\F$ has a ring language on it along with the names for elements of the subfield $\kk.$ In particular, $\X$ is a imaginary sort in $\F_{|\kk}$ (a definable subset of $\PP^N(\F)$).  $\mu_\alpha:\X_{\mu,\alpha}\to \X$ are definable using parameters $\alpha.$

The group $ \mathrm{GDeck}(\X_{\mu,\alpha}/\X)$ is another structure interpretable in  $\F_{|\kk}$ using parameter $\alpha,$ uniformly on  $\alpha\in \mathrm{Zeros}\,\mathbf{f}_{\mu}.$


\item  For each $\mu\in \mathcal{M}_\X,$ $\alpha\in \mathrm{Zeros}(\mathbf{f}_\mu)$ there is a binary function symbol $*$ defining the action 
$$ \mathrm{GDeck}(\X_{\mu,\alpha}/\X)\times \tilde{\pp}_{\mu,\alpha}\to \tilde{\pp}_{\mu,\alpha},\ \ 
(g,p)\mapsto g*p.$$

\item 
For each $\mu\in \mathcal{M}_\X,$ there is a relation symbol $\Phi_\mu$ between $\U,$ $\F$ and $\tilde{\pp}_\mu$ which, for $\alpha\in  \mathrm{Zeros}(\mathbf{f}_\mu),$   
 defines a map
$$\Phi_{\mu,\alpha}: \tilde{\pp}_{\mu,\alpha}\times \U\to \X_{\mu,\alpha}.$$

Given $p\in \tilde{\pp}_{\mu,\alpha},$ which has an interpretation of  a map
$\U\to \X_{\mu,\alpha}$ we set
$$\Phi_{\mu,\alpha}(p,u)=p(u).$$

\end{itemize}


\medskip

It is clear from the definition of $\mathcal{L}_\X$ that any model $(\U,\F_{/\kk})$  of $T^\sharp_\X$ can be transformed into a structure in the language $\mathcal{L}_\X$ by adding certain 0-definable maps and
sorts of 
$(\U,\F_{/\kk}, \{ \tilde{\pp}_\mu\}_{\mu\in \MM})$ and forgetting the names for  maps $\pp_{\mu,\alpha}$ of $\mathcal{L}^\sharp_\X.$

\medskip

We denote $\XX^{et}(\F)$ the structure in the language $\mathcal{L}_\X$ which corresponds in this way to $\XX^{et}_\sharp(\F)$ considered in \ref{projlimit}.

\epk

\bpk \label{Tflat} Define the {\bf theory $T_\X$ in the language $\mathcal{L}_\X$} by the following axioms:
\begin{itemize}
\item[A1] {\em $\F$ is an algebraically closed field of characteristic $0$ with subfield $\kk$ of constants.}


\item[A$2'(\mu$)]  For the given $\mu\in \mathcal{M}_\X,$ for any  zero $\alpha$ of $\mathbf{f}_\mu$ and $p\in \tilde{\pp}_{\mu,\alpha},$ the map  
 $$u\mapsto \Phi_{\mu,\alpha}(p,u), \ \ \U\to \X_{\mu,\alpha}$$ 
is surjective.

For any $g\in \mathrm{GDeck}(\X_{\mu,\alpha}/\X)$
$$ \Phi_{\mu,\alpha}(g*p,u)=g( \Phi_{\mu,\alpha}(p,u)).$$

 For any $p_1,p_2\in \tilde{\pp}_{\mu,\alpha}$ there is $g\in \mathrm{GDeck}(\X_{\mu,\alpha}/\X)$ such that $g*p_1=p_2.$ 
In particular, $\tilde{\pp}$ is a set with single element $\pp.$ 

\item[A3'] For given $\mu\in \mathcal{M}_\X$ and a zero $\alpha$ of $\mathbf{f}_\mu$ 
there is $q\in \tilde{\pp}_{\mu,\alpha}$ such that
$$\pp\circ q\inv =\mu_\alpha.$$

\end{itemize}  
\epk
\bpk \label{reducts} {\bf Lemma.} {\em Given a model $(\U, \F_{/\kk}, \{ \pp_\mu\}_{\mu\in \MM})$ of $T^\sharp_\X$
there is a model $(\U, \F_{/\kk}, \{ \tilde{\pp}_\mu\}_{\mu\in \MM})$ of  $T_\X$ interpretable in $(\U, \F_{/\kk}, \{ \pp_\mu\}_{\mu\in \MM}).$
}

{\bf Proof.} The interpretation is just by definition \ref{ppp}:\  the set $\tilde{\pp}_{\mu,\alpha}$ is in bijective correspondence with the set $\mathrm{GDeck}(\X_{\mu,\alpha}/\X)\times \{ \pp_{\mu,\alpha}\}$ and $\tilde{\pp}_{\mu}$ with the family $\mathrm{GDeck}(\X_{\mu,\alpha}/\X)\times \{ \pp_{\mu,\alpha}\}:\ \mathbf{f}_\mu(\alpha)=0.$ 
$\Box$
\epk
\bpk {\bf Theorem.} {\em $T_\X$ is a complete theory allowing elimination of quantifiers.

}

{\bf Proof.} Any model of   $T_\X$ can be made into a model of $T^\sharp_\X$ by
setting 
\be \label{choice} \pp_{\mu,\alpha}(u):= \Phi_{\mu,\alpha}(q,u)\ee
for some choice of $q\in \tilde{\pp}_{\mu,\alpha},$ which is possible by axiom A3'.

It follows that   any two saturated model of   $T_\X$ of the same cardinality are isomorphic and hence the completeness. 

Elimination of quantifiers follows by the same back-and-forth  construction in the proof of \ref{complete} in the language $\mathcal{L}_\X.$
$\Box$

\epk

From now on {\bf we work in the language $\mathcal{L}_\X.$}

For each $\mu\in \mathcal{M}_\X,\ \alpha\in \mathrm{Zeros}\, \mathbf{f}_\mu$ we fix $b_{\mu,\alpha}\in \kk_\alpha(\X_{\mu,\alpha})$ which generates the function field over $\K(\alpha),$ that is 
$\K_{\mu,\alpha}=\K(\alpha,b_{\mu,\alpha}).$ 
We set $\tilde{b}_{\mu,\alpha}$ to be the orbit of $b_{\mu,\alpha}$ under the Galois group $\mathrm{Gal}(\K(\alpha,b):\K(\alpha)).$ We may identify $b_{\mu,\alpha}$ and its conjugates as
 rational functions 
 $b: \X_{\mu,\alpha}(\F)\to \F$
defined over $\kk(\alpha),$ with domain of definition dense in the variety.

Note that applying $\sigma \in \mathrm{Gal}(\kk(\alpha):\kk)$ 
to $b_{\mu,\alpha}$ and to  $\X_{\mu,\alpha}$ we obtain a rational function
$b_{\mu,\alpha'}: \X_{\mu,\alpha'}(\F)\to \F$ where $\alpha'=\sigma(\alpha)\in \kk(\alpha),$ and so
$b_{\mu,\alpha'}\in  \kk({\alpha'})(\X_{\mu,\alpha'}).$

 \bpk\label{ia} {\bf Lemma.} {\em There is a  
 family $i_{\mu,\alpha}: \mu\in \mathcal{M}_\X,\ \alpha\in \mathrm{Zeros}\, \mathbf{f}_\mu$ of bijections
$$i_{\mu,\alpha}:\tilde{\pp}_{\mu_\alpha}\rightarrow \tilde{b}_{\mu_\alpha}\subset \K_{\mu,\alpha}$$ which is an isomorphism between the finite structures induced by the ambient structures in the category $\mathfrak{M}$ (a bi-interpretation) over $\kk(\alpha).$ 

Moreover, there is a $L_{\omega_1,\omega}$-interpretation of a field $\mathrm{P}_{\mu,\alpha}$  in structure
 $\tilde{\pp}_{\mu_\alpha}$ so that
$\tilde{\pp}_{\mu_\alpha}\subset \mathrm{P}_{\mu,\alpha}$
and $i_{\mu,\alpha}$ can be extended to an isomorphism of fields
$$
i_{\mu,\alpha}: \mathrm{P}_{\mu,\alpha}\to \K_{\mu,\alpha}.$$

}

{\bf Proof.} 
Note that  (\ref{GDeck}) of \ref{defDeck} asserts  the existence of an isomorphism  
$$j_{\mu,\alpha}: \mathrm{GDeck}(\X_{\mu,\alpha}/\X)\to
\mathrm{Gal}(\K(\alpha,b_{\mu,\alpha}):\K(\alpha)).$$

Clearly, given $g\in \mathrm{GDeck}(\X_{\mu,\alpha}/\X)$ we get a Galois automorphism of the field of rational functions $\K_{\mu,\alpha}$ over $\kk_\alpha,$ 
$$\hat{g}: b\mapsto b\circ g.$$
Thus we may assume $j_{\mu,\alpha}(g)=\hat{g}.$
  
Set, for $b\in \tilde{b}_{\mu,\alpha},$
 $$i_{b}: g\circ \pp_{\mu,\alpha}\mapsto 
b\circ g, \ \ \ g\in \mathrm{GDeck}(\X_{\mu,\alpha}/\X).$$   

This is injective and gives us
$$\tilde{b}_{\mu,\alpha}=i_{b}(\tilde{\pp}_{\mu,\alpha}).$$

Set $i_{\mu,\alpha}:= i_{b_{\mu,\alpha}}.$

The $\kk(\alpha)$-definable   relations between elements $b_1,\ldots,b_k$ of $\tilde{b}_{\mu,\alpha}$ induced from the ambient field  (equivalently, the relation invariant under $\mathrm{Gal}(\K(\alpha,b_{\mu,\alpha}):\K(\alpha))$ are boolean combinations of relations of the form 
$f(\alpha, b_1,\ldots,b_k)=0,$ where $f$ is a polynomial over $\kk.$

Set a relation between $p_1,\ldots,p_k\in \tilde{\pp}_{\mu,\alpha}$ and $\alpha$
$$R_f(\alpha,p_1,\ldots,p_k):\equiv  f(\alpha,i_{b_{\mu,\alpha}}(p_1),\ldots,i_{b_{\mu,\alpha}}(p_k))=0.$$

Note that 
$$f(\alpha,i_{b_{\mu,\alpha}}(p_1),\ldots,i_{b_{\mu,\alpha}}(p_k))=0 \Leftrightarrow \exists b\in \tilde{b}_{\mu,\alpha}\ f(\alpha,i_{b}(p_1),\ldots,i_{b}(p_k))=0$$
since $i_{b_{\mu,\alpha}}(p)\mapsto i_b(p)$ is a $\mathrm{Gal}(\K(\alpha,b_{\mu,\alpha}):\K(\alpha))$-transformation.

Thus, $R_f$ is 0-definable (definable over $\kk$). 

Now we interpret the field structure $\mathrm{P}_{\mu,\alpha}$ in the substructure with the universe $\tilde{\pp}_{\mu,\alpha}$ using language $L_{\omega_1,\omega}$ as follows:

The universe of  $\mathrm{P}_{\mu,\alpha}$ will be interpreted as $S_{\mu,\alpha}/E_{\mu,\alpha}$ where $S_{\mu,\alpha}$
 is the $L_{\omega_1,\omega}$-definable set consisting of formal terms
$F(\alpha,p_1,\ldots,p_N),$ for $F\in \kk[X_0,X_1,\ldots, X_N],$ $p_1,\ldots,p_N$ the list of all elements of $\tilde{\pp}_{\mu,\alpha},$ and $E_{\mu,\alpha}$ is the equivalence relation between the terms $F_1,F_2,$
$$R_{F_1-F_2}(\alpha,p_1,\ldots,p_N).$$
That is 
$ F_1(\alpha,p_1,\ldots,p_N)-F_2(\alpha,p_1,\ldots,p_N)=0$ when interpreted by $i_{\mu,\alpha}.$ 

The operations $+$ and $\times$ on the set of terms 
$S_{\mu,\alpha}$ gives it the structure of a ring. And taking the quotient by $E_{\mu,\alpha}$ we get, by construction, field $\mathrm{P}_{\mu,\alpha}$ isomorphic to $\K_{\mu,\alpha}.$ Clearly,
$\tilde{\pp}_{\mu,\alpha}\subset \mathrm{P}_{\mu,\alpha}$ since the equivalence $E_{\mu,\alpha}$ is trivial on $\tilde{\pp}_{\mu,\alpha}.$

Finally, we claim that  $\tilde{\pp}_{\mu,\alpha}$ with the structure induced from $\XX^{et}(\F)$ is interpreted in the field structure  $\mathrm{P}_{\mu,\alpha}.$  We first note that by the construction of $\XX^{et}(\F)$ in \ref{ppp} - \ref{Lflat} the transformation of $\tilde{\pp}_{\mu,\alpha},$ $p\mapsto g\circ p,$ by the action of a $g\in \mathrm{GDeck}(\X_{\mu,\alpha}/\X),$
can be extended to an automorphism of $\XX^{et}(\F)$ fixing $\kk(\alpha).$ This implies that definable relations on $\tilde{\pp}_{\mu,\alpha}$ are invariant under the action. Equivalently, the image of such relation under $i_{\mu,\alpha}$ is invariant under the action by the Galois group of the function field. Thus definable relations are boolean combinations of the $R_f,$ which proves the claim. 
 
$\Box$

\epk  

\bpk \label{varphi} {\bf Lemma.} {\em Let  
$$\mu_\alpha: \X_{\mu,\alpha}\to \X,\ \ \nu_\beta: \X_{\nu,\beta}\to \X,\ \ \mu, \nu\in \mathcal{M}^0_\X, $$
with a matching covering
$$(\mu_\alpha^{-1}\nu_\beta): \X_{\nu,\beta}\to  \X_{\mu,\alpha}.$$

Let $$(\mu_\alpha^{-1}\nu_\beta)^*: \K_{\mu,\alpha} \to \K_{\nu,\beta} $$
the embedding of fields induced by covering morphism $(\mu_\alpha^{-1}\nu_\beta).$

Let   $$i_{\mu,\alpha}:\mathrm{P}_{\mu,\alpha}\rightarrow  \K_{\mu,\alpha}$$
as constructed in \ref{ia}. 

Then one can adjust the construction of 
 $$i_{\nu,\beta}:\mathrm{P}_{\nu_\beta}\rightarrow \K_{\nu,\beta}$$
so that the diagram commutes}


\begin{center}

 \begin{tikzpicture}
  
 \draw (10,8) node {$i_{\mu,\alpha}$};
\draw (8,7.5) node {$\mathrm{P}_{\mu,\alpha}$};
\draw[->][line width=0.3mm] (8,7 )--(8, 5.3);\draw[->][line width=0.3mm] (12.1,7 )--(12.1, 5.3);
 \draw (7.2,6.4) node {$(\mu_\alpha^{-1}\nu_\beta)^*$};  \draw (13,6.4) node {$(\mu_\alpha^{-1}\nu_\beta)^*$};
\draw[->][line width=0.3mm] (8.5, 7.5)--(11.2, 7.5);
\draw (12.2,7.5) node {$\K_{\mu,\alpha}$};
\draw (12.2,5) node {$\K_{\nu,\beta}$};

 \draw (10,5.3) node {$i_{\nu,\beta}$};
\draw (8,5) node {$\mathrm{P}_{\nu,\beta}$};

\draw[->][line width=0.3mm] (8.5, 5)--(11.2, 5);




\end{tikzpicture}

\end{center} 

{\bf Proof.} Let $b_{\mu,\alpha}$ and $b_{\nu,\beta}$ be the generating elements of function fields as above. The embedding  $\K_{\mu,\alpha}\subseteq \K_{\nu,\beta}$
gives rise to a $\kk$-rational function $h=\frac{F_1}{F_2},$ $F_1,F_2\in \kk[X_0,X_1],$ such that
$h(\beta,b_{\nu,\beta})=b_{\mu,\alpha}$ and so $b\mapsto h(\beta,b)$ is the restriction of the map $(\mu_\alpha^{-1}\nu_\beta)^*,$  $\tilde{b}_{\nu,\beta}\to \tilde{b}_{\mu,\alpha},$
which extends to embedding of respective fields.

Then in the field $\mathrm{P}_{\nu,\beta}$ 
 of formal terms embeds into   $\mathrm{P}_{\mu,\alpha}$ using the same rational function
 $h=\frac{F_1}{F_2}.$ $\Box$

\epk

\bpk \label{tildeK} Since the family of fields $$\mathcal{M}_\K:=\{ \K_{\mu,\alpha}: \mu\in \mathcal{M}_\X, \alpha\in \mathrm{Zeros}\,\mathbf{f}_\mu\}$$
forms an inverse system with reprect to embeddings 
$(\mu_\alpha^{-1}\nu_\beta)^*,$ the following inverse lmit is well-defined 

$$\tilde{\K}= \lim_{\leftarrow\mathcal{M}}\K_{\mu,\alpha} $$ 

By this definition   $\tilde{\K}$ is the union of all the function fields of Galois \'etale covers of $\X.$ When speaking of it as a structure we consider it a field over $\K,$ that is with elements of $\K$ in $\tilde{\K}$ being names. Note that it automatically names also elements of $\kk.$

\medskip

We give names and consider the multisorted structures definable or interpretable in $\XX^{et}(\F):$
$$\tilde{\pp}_\X:=\{ \tilde{\pp}_{\mu}: \mu \in \mathcal{M}\},\ \ \mathbf{P}_\X:= \{ \mathbf{P}_\mu: \ \mu\in \mathcal{M}\},\ \ \mathbf{P}_\mu:= \{ \mathbf{P}_{\mu,\alpha}: \alpha\in \mathrm{Zeros}\, \mathbf{f}_\mu\}$$
with relations induced from the ambient structure.

\medskip

Set $\XX^{et}_{fin}(\kk^{alg})$ to be the substructure of $\XX^{et}(\kk^{alg})$ obtained by removing the covering sort $\tilde{\X}^{et}(\kk^{alg})$ but keeping all the induced relations.

 \epk
\bpk \label{M}   {\bf Theorem.} {\em 

1.The field $\tilde{\K}$ as a structure  is isomorphic to
$\mathbf{P}_\X$ via an isomorphism
$$\mathbf{P}_\X\cong_{i_\X} \tilde{\K}.$$

2. $\mathbf{P}_\X$ is bi-interpretable with $\tilde{\pp}_\X,$ i.e. $$\mathbf{P}_\X
 \cong_\mathfrak{M} \tilde{\pp}_\X.$$

3.
  $$\Aut(\mathbf{P}_\X)
 \cong \Aut( \XX^{et}(\kk^{alg}))$$
and the isomorphism is induced by the bi-interpretation 
$$\mathbf{P}_\X\cong_\mathfrak{M}  \XX^{et}_{fin}(\kk^{alg}).$$

 }

{\bf Proofs.} 1. Follows  directly from \ref{varphi}. 

\medskip

2. Follows directly from \ref{ia}.

\medskip

3. First, we prove that
$$\dcl (\mathbf{P}_\X)\supset \X_{\mu,\alpha}(\kk^{alg})\mbox{ for all }\mu,\alpha.$$
where $\dcl$ is in the first-order language.

Indeed, note that $$\kk^{alg}=\bigcup\{ \kk(\alpha): \mathbf{f}_\mu(\alpha)=0, \mu\in \mathcal{M}\}$$
and each such $\alpha$ is definable as a code for the sort $\tilde{\pp}_{\mu,\alpha}$ in $\mathbf{P}_\X.$ Thus $\dcl (\mathbf{P}_\X)\supset \kk^{alg}$ and so $\dcl (\mathbf{P}_\X)\supset \X_{\mu,\alpha}(\kk^{alg})$ for all $\mu,\alpha.$

Now we conclude that  $$\dcl (\mathbf{P}_\X)\supseteq  \XX^{et}_{fin}(\kk^{alg})$$ since
the right-hand side consists of sorts $\X_{\mu,\alpha}(\kk^{alg})$ and $\tilde{\pp}_{\mu,\alpha}$ only.  In particular, any $\sigma\in \Aut(\mathbf{P}_\X)$ induces a unique automorphism of $\XX^{et}_{fin}(\kk^{alg}).$ 

 Finally, we claim that, in its own turn,  any  $\sigma\in  \Aut( \XX^{et}_{fin}(\kk^{alg}))$ induces a unique automorphism of $ \XX^{et}(\kk^{alg})),$ that is has a unique extension to   $\tilde{\X}^{et}(\kk^{alg})),$

Indeed, consider $u\in \tilde{\X}(\kk^{alg}).$ Choose
elements $\pp_{\mu,\alpha}\in \tilde{\pp}_{\mu,\alpha}$ and set $x_{\mu,\alpha}:=\pp_{\mu,\alpha}(u),$  element of $\X_{\mu,\alpha}(\kk^{alg}).$ By construction $u$ is the unique element satisfying the system of equations $$x_{\mu,\alpha}=\pp_{\mu,\alpha}(u):\ \mu\in \mathcal{M}, \mathbf{f}_\mu(\alpha)=0,$$
i.e. $u$ is type definable over   $\XX^{et}_{fin}(\kk^{alg})).$ 
Hence, the action of $\sigma$ on $\XX^{et}_{fin}(\kk^{alg}))$  has a unique extension to   $\tilde{\X}^{et}(\kk^{alg})).$


Note that $\mathbf{P}_\X\hookrightarrow \XX^{et}(\kk^{alg})$
by definition. Together with the claim this completes the proof. $\Box$

\epk

\bpk\label{Main2} {\bf Corollary.} 
$$\Aut\, \XX^{et}(\kk^{alg})\cong\Aut\, \XX^{et}_{fin}(\kk^{alg})\cong \Aut\, \tilde{\pp}_\X\cong   \mathrm{Gal}(\tilde{\K}/\K).$$

Note that  $\XX^{et}_{fin}(\kk^{alg})$ is a finitary structure.
\epk

\bpk\label{fact} {\bf Fact} (\cite{Lenstra}, Corollary 6.17).
{\em For a normal scheme $\X$ with a function field $\K$ 
$$\pi_1^{et}(\X,x)=\mathrm{Gal}(\tilde{\K}/\K),$$
where $\pi_1^{et}(\X,x)$ is the \'etale fundamental group as defined in} \cite{SGA1}.


\epk

\bpk {\bf Corollary.}
$$\pi_1^{et}(\X,x)\cong  \Aut(\XX^{et}_{fin}(\kk^{alg})).$$

\epk
\bpk {\bf Remark.}
The embedding $\F_{/\kk}\hookrightarrow \XX^{et}(\F)$ induces via the functor $\Aut$ of section \ref{s1} and Proposition \ref{Pr1}
$$\hat{\Gamma}\hookrightarrow \Aut(\XX^{et}(\F))\twoheadrightarrow \Aut \F_{/\kk}$$
in particular, when $\F=\kk^{alg},$
$$\hat{\Gamma}\hookrightarrow \pi_1^{et}(\X,x)\twoheadrightarrow \mathrm{Gal}_\kk.
$$

\epk
The following gives a link with a general model theory setting. 

\bpk {\bf Proposition.} {\em  
$\pi_1^{et}(\X,x)$ is isomorphic to the {\bf Lascar group} of theory $T_\X.$}

{\bf Proof.} Lascar group for stable theories  is known to be isomorphic to $\Aut(\acl^{eq}(0))$ of a model. Hence 
it is enough to prove that 
$\acl^{eq}(0)$ is bi-interpretable with
the substructure $\tilde{\pp}_\X$ of any model.
Equivalently, $\acl^{eq}(0)=\dcl (\tilde{\pp}_\X).$ 

The inclusion 
$\acl^{eq}(0)\supseteq \tilde{\pp}_\X$ is obvious since $\tilde{\pp}_\X$ is finitary (the union of finite sorts). To prove the inverse we can use the language which names all elements $\pp_{\mu,\alpha}$ in sorts  $\tilde{\pp}_{\mu}.$ This language is equivalent to $\mathcal{L}_\X^\sharp.$ and so we can use  theorem \ref{complete} describing definable sets in models of the theory. It is easy to see that the only finite imaginary sorts are the ones on finite sorts $\tilde{\pp}_{\mu}$ and, by elimination of imaginaries in algebraically closed fields, subsets of $\kk^{alg}.$ In terms of language $\mathcal{L}_\X$ both are part of  
$\acl^{eq}(0).$
 $\Box$  
\epk
\bpk {\bf Definition of $\pi_1^{et}(\X,x).$}

Let $x\in \X(\F)$ and consider
 the  multi-sorted structure 
 $$\mathrm{F}_{x}=\{ \mu_\alpha\inv(x): \mu\in \mathcal{M},\ \mathbf{f}_\mu(\alpha)=0\}$$
 with relations induced on it from $\XX^{et}(\F).$ 

Claim.  There is an $x$-definable bijection $i_{x}: \tilde{\pp}_\X\to \mathrm{F}_{x}$ which induces an interpretation of $ \tilde{\pp}_\X$ in  $\mathrm{F}_{x}.$
In particular,   $$i_{x}(\tilde{\pp}_{\mu,\alpha})=\mu_\alpha\inv(x).$$

Proof. Note that for each $\mu_\alpha,$ 
$$ \mu_\alpha\inv(x)=\pp_{\mu,\alpha}(\pp\inv(x))=p(\pp\inv(x))\mbox{ for any }p\in  \tilde{\pp}_{\mu,\alpha}.$$
Note that $ \pp\inv(x)$ is $L_{\omega_1,\omega}$-definable in  $\mathrm{F}_{x}$
as the projective limit of  $\{ \mu_\alpha\inv(x)\}.$

For each $\mu_\alpha,$ 
$$ \mu_\alpha\inv(x)=\pp_{\mu,\alpha}(\pp\inv(x))=p(\pp\inv(x))\mbox{ for any }p\in  \tilde{\pp}_{\mu,\alpha}.$$
Thus one can set, for each $p\in \tilde{\pp}_{\mu,\alpha},$
 $$i_{x}(p):  \pp\inv(x)\to  \mu_\alpha\inv(x),  \ u\mapsto p(u),$$ the restrictions of respective functions $\tilde{\X}(\F)\to \X_{\mu,\alpha}(\F).$

The relations between the  $\tilde{\pp}_{\mu,\alpha}$ descend to relations on $\mathrm{F}_{x}$ in the obvious way. Claim proved.

Define $\mathrm{F}_{x}^\mathrm{forget}$ to be the reduct of $\mathrm{F}_{x}$ to the language expressing the relations on $\tilde{\pp}_\X$ only.  Then
$$\tilde{\pp}_\X\cong_\mathfrak{M} \mathrm{F}_{x}^\mathrm{forget}$$
via the bijection $i_{x}.$

It follows that $\Aut\,  \tilde{\pp}_\X\cong \Aut\,\mathrm{F}_{x}^\mathrm{forget}$ via the bijection $i_{x}.$ Define
$$\pi_1^{et}(\X,x)=\Aut\,\mathrm{F}_{x}^\mathrm{forget}.$$
 \epk
 
Finally we use our construction to reformulate Grothendieck's section conjecture.

\bpk \label{8.1} {\bf Theorem.} {\em Suppose $\kk^\sharp=\kk.$ Then there is a section 
$$s: \mathrm{Gal}_\kk\to \pi_1^{et}(\X)$$
of the canonical homomorphism  $ \pi_1^{et}(\X)\to \mathrm{Gal}_\kk.$
  }

{\bf Proof.}

By assumption  $$\Aut \XX^{et}_\sharp(\kk^{alg})=\mathrm{Gal}_\kk.$$

We can interpret  $\XX^{et}(\kk^{alg})$ 
in   $\XX^{et}_\sharp (\kk^{alg})$ without parameters as explained in \ref{ppp} and \ref{Lflat}. It follows that any automorphism of 
 $\XX^{et}_\sharp (\kk^{alg})$ extends uniquely to an  automorphism of 
 $\XX^{et} (\kk^{alg}).$ That is there is an embedding 
 $$s: \Aut\, \XX^{et}_\sharp (\kk^{alg})\to \Aut\, \XX^{et} (\kk^{alg}).$$
 Under our assumption this is an embedding 
  $$s: \mathrm{Gal}_\kk\to \Aut\, \XX^{et} (\kk^{alg})$$
  which is obviously a section of the restriction homomorphism
$$ \Aut\, \XX^{et} (\kk^{alg})  \to  \mathrm{Gal}_\kk.$$ $\Box$

\epk
Now we formulate a version of section conjecture which is formally stronger than the conventional section conjecture.
\bpk \label{8.3} {\bf $\sharp$-version of section conjecture.} Suppose there is a section
$$s: \mathrm{Gal}_{\kk}\to \Aut\, \XX^{et}_\sharp(\kk^{alg}).$$
Then $\X$ has a $\kk$-rational point and $\kk^\sharp=\kk.$
\epk

\include{Plim}

\end{document}